\documentclass[11pt]{article}

\usepackage{amssymb}
\usepackage{amsmath}
\usepackage{graphics}

\newcommand{\qed}{$\;\;\;\Box$}
\newenvironment{proof}{\par\smallbreak{\sl\bf Proof.~}}
{\unskip\nobreak\hfill \qed \par\medbreak}

\newcounter{claim}
\renewcommand{\theclaim}{\arabic{claim}}
\newenvironment{claim}{\refstepcounter{claim}%
\par\medskip\par\noindent{\bf Claim~\theclaim.}\rm}%
{\par\medskip\par}

\newenvironment{subproof}{\par\noindent{\bf Proof of Claim.}}%
{\qed\par\smallbreak}
\renewcommand{\L}{{\mbox{L}}}

\newcommand{\D}{{\cal D}}
\newcommand{\E}{{\cal E}}

\newcommand{\N}{{\mathbb N}}
\newcommand{\R}{{\mathbb R}}

\newcommand{\beq}{\begin{equation}}
\newcommand{\ee}{\end{equation}}

\renewcommand{\d}{\partial}

\newtheorem{theorem}{Theorem}
\newtheorem{lemma}[theorem]{Lemma}

\newtheorem{defn}[theorem]{Definition}
\newtheorem{corollary}[theorem]{Corollary}
\newtheorem{remark}[theorem]{Remark}

\newcommand{\al}{\alpha}
\newcommand{\be}{\beta}
\newcommand{\ga}{\gamma}

\newcommand{\eps}{\varepsilon}
\newcommand{\vphi}{\varphi}
\newcommand{\la}{\lambda}
\newcommand{\om}{\omega}
\newcommand{\io}{\iota}

\newcommand{\e}{\mathop{\rm exp}}

\date{}

\title{
Smoothing Solutions to Initial-Boundary Problems for 
First-Order Hyperbolic Systems} 

\author{
Irina Kmit\thanks{Institute of Mathematics, Humboldt University of Berlin,
Rudower Chaussee 25,
D-12489 Berlin, Germany and Institute for Applied Problems of Mechanics and Mathematics,
Ukrainian Academy of Sciences,  Naukova St.\ 3b,
79060 Lviv,
Ukraine.
{\small   
E-mail:
{\tt kmit@informatik.hu-berlin.de}}
}
}

\begin{document}

\maketitle

\begin{abstract}
We consider initial-boundary problems for general linear first-order strictly
hyperbolic systems
with local or nonlocal nonlinear boundary conditions. While boundary data 
are supposed to be smooth, initial conditions can contain distributions of any
order of singularity. It is known that such problems have a unique continuous solution
if the initial data are continuous. 
In the case of strongly singular initial data we prove the existence  of a 
(unique) delta wave solution. In both cases, we say that a solution is smoothing if 
it eventually becomes $k$-times continuously differentiable for each $k$.
Our main result is a criterion allowing us to determine
whether or not the solution is smoothing. In particular, 
we prove a rather general smoothingness result in
the case of classical boundary conditions.


\end{abstract}

\section{Introduction}\label{sec:intr}
Solutions to  hyperbolic PDEs demonstrate a wide spectrum of regularity
behavior. The appearance of singularities in nonlinear cases is known as the blow up of a 
solution~\cite{alinhac,hoermander}. Singularities can appear
in a finite time even for small and smooth initial data~\cite{John}. In some cases, both linear and 
nonlinear, a solution with time either preserves  the same regularity as it has on the 
boundary
or becomes less or more
regular in time. The singularities encountered in the latter case are called anomalous~\cite{popivanov06,popivanov03}. 
Criteria for the appearance of anomalous singularities are given in~\cite{Elt,LavLyu,Ober86,popivanov03,RauchReed81}.
These papers are devoted to the case of interaction between singularities weaker than 
 the Dirac measure. The singularities resulting from this interaction turn out to be weaker than the
incoming singularities. A different effect is observed in~\cite{distr,travers}, where
the incoming singularities are derivatives of the Dirac measure. In this case the  interaction 
produces singularities stronger than the initial ones. We will focus on the phenomenon 
of improving regularity in the case of initial-boundary 
value problems  with nonlinear local and
nonlocal boundary conditions for first-order linear strictly hyperbolic systems.

Specifically, in the domain $\Pi=\{(x,t)\,|\,0<x<1, t>0\}$
we address the problem
\begin{equation}\label{eq:1}
(\partial_t  + \Lambda(x,t)\partial_x + A(x,t))
 u = g(x,t), \qquad\ (x,t)\in \Pi
\end{equation}

\begin{eqnarray}
&&u(x,0) = \varphi(x),\quad\qquad\ \, x\in (0,1)\label{eq:2}\\&&
\begin{array}{ll}
u_i(1,t) = h_i(t,v(t)), & 1\le i\le k, \ \  t\in(0,\infty)\\
u_i(0,t) = h_i(t,v(t)), & k< i\le n, \ \  t\in(0,\infty),
\end{array}\label{eq:3}
\end{eqnarray}
where  $u$, $g$, and $\varphi$  are real $n$-vectors, $A=\{a_{ij}\}_{i,j=1}^n$,
$\Lambda=diag(\lambda_1,\dots,\lambda_n)$, and
\begin{equation}\label{eq:vv}
v(t)=\left(v_1(t),\dots,v_n(t)\right)=\left(u_1(0,t),\dots,u_{k}(0,t),u_{k+1}(1,t),\dots,u_{n}(1,t)\right).
\end{equation}
Note that boundary conditions~(\ref{eq:3})  cover the cases of classical boundary
conditions (if $h_i$ do not depend on $v$) and reflection boundary conditions of local
and nonlocal type.

We assume that
\begin{equation}\label{eq:L1}
\lambda_1<\dots<\lambda_k<0<\lambda_{k+1}<\dots<\lambda_n
\end{equation}
for all $(x,t)\in\overline\Pi$. Condition  (\ref{eq:L1}) occurs
in many applications, where the functions $u_j$ for $j\le k$ (resp. $k+1\le j\le n$)
 describe the ``species'' that travel to the left (resp. to the right) 
and are reflected in $x=0$ (resp. $x=1$) according to 
the boundary conditions (\ref{eq:3}). 

We will impose the following smoothness assumptions on the initial data: The entries of $\Lambda$,
$A$, $g$, and $h=(h_1,\dots,h_n)$ are smooth in all their arguments in the respective domains,
while the entries of $\vphi$ are allowed to be either continuous functions or strongly 
singular distributions.

In the case of a continuous $\vphi$ (considered in Section~2),
by a solution to problem (\ref{eq:1})--(\ref{eq:3}) we mean a 
{\it continuous solution}, i.e., a continuous vector-function in $\overline\Pi$
 satisfying an integral system equivalent to
(\ref{eq:1})--(\ref{eq:3}). The existence and uniqueness of a continuous
solution is proved in~\cite{ijdsde} (see Theorem~\ref{thm:ijdsde}).
In the case of a strongly singular $\vphi$ 
(considered in Section~4) by a solution we mean a {\it delta wave
solution}, i.e., a weak limit of  solutions to the
original problem with regularized  initial data that does
not depend on a particular regularization. We refer the reader 
to~\cite{Obe,RauchReed81} for a more detailed definition and motivation of delta waves. 
Theorem~\ref{thm:delta} in Section~4 establishes the existence of a delta wave
solution for a version of (\ref{eq:1})--(\ref{eq:3}).

It is clear that the regularity of initial conditions
(\ref{eq:2}) constraints the regularity of a solution
if the latter is considered
in the entire domain $\Pi$.  However, the influence of the initial data can be suppressed
if the regularity behavior is considered in dynamics, starting from a point of time $T$. 

\begin{defn}\rm
A solution $u$ to problem (\ref{eq:1})--(\ref{eq:3}) is called {\it smoothing} 
if, whatever 
$m\in\N$, there exists $T>0$ such that
$u\in C^m\left(\overline\Pi\cap\{t\ge T\}\right)^n$.
\end{defn}

Our main result is a smoothingness criterion for solutions to
problem (\ref{eq:1})--(\ref{eq:3})
in terms of the layout of characteristic curves 
(Theorems~\ref{thm:I} and~\ref{thm:reg-delta}).
In the case of classical boundary condition, the criterion implies that the solution is smoothing
whenever $\inf|\lambda_i|>0$ for all $i\le n$. 
As another consequence, we obtain 
 a class of boundary conditions under which the wave equation has smoothing solutions
(see~\cite{wave} for a special case of this result).

Our analysis of problem (\ref{eq:1})--(\ref{eq:3}) shows a phenomenon usually
observed in the situations when solutions to hyperbolic PDEs change their 
regularity: the smoothness changes jump-like rather than gradually.
Another feature of  problem (\ref{eq:1})--(\ref{eq:3}) shown in~\cite{ijdsde}
is that, even if we allow non-Lipschitz nonlinearities in (\ref{eq:3}), 
this system demonstrates almost linear behavior. The smoothingness effect in the  non-Lipschitz setting 
contrasts with blowups~\cite{alinhac,hoermander} observed in many nonlinear systems.

In \cite{Elt,LavLyu}, results similar to ours are obtained in some more special cases, namely
for  homogeneous linear boundary conditions of local type with constant coefficients and 
continuous initial data. Some restrictions on
(\ref{eq:1})--(\ref{eq:3}) are imposed in~\cite{Elt,LavLyu} by technical reasons, as the authors use
an approach based on
the Laplace transformation and the Green's function method.
We extend these results  to the case of general linear first-order hyperbolic systems, 
nonlinear nonlocal boundary conditions, and distributional 
initial data of the Dirac delta type and derivatives thereof.
Note that we use a different approach based on the classical method of
characteristics. This method suits well for understanding of the
mechanism of the smoothingness effect.

An essential technical difficulty to overcome in 
demonstrating the regularity self-improvement is caused by the fact that  
the domain of influence of initial conditions (\ref{eq:2})  is in general 
infinite. In other words, the regularity of solutions all the time depends on the
regularity of the initial data. However, this dependence is different for the
boundary and the integral parts of the equivalent integral form 
of  problem (\ref{eq:1})--(\ref{eq:3}). Since the boundary
summands are compositions of boundary data with functions defining characteristic curves
and hence are ``responsible'' for propagation of singularities, the smoothingness effect is encountered whenever  
the boundary summands have a ``bounded memory'' or, more rigorously, 
all characteristics of (\ref{eq:1}) are bounded and
each boundary singularity expands inside $\Pi$ 
along a finite number of characteristic curves.\footnote{This sharply  contrasts with the case of a 
Cauchy problem where  solutions cannot be 
smoothing because the boundary part all the time "remembers" 
the regularity of the initial data.} Our strategy of obtaining the
smoothingness criterion consists in identifying an appropriate class of
boundary conditions ensuring the "bounded memory" property and in showing
that  the integral summands do eventually improve upon the regularity of the
initial data. It is important for our analysis that the lower order terms,
with the exception of the diagonal ones,
 contribute into the integral summands transversely to
characteristic directions. This is ensured by (\ref{eq:L1}). Therefore, the integral part of the system 
 causes no propagation of singularities
and, moreover, suppresses it.

The mathematical motivation of the paper is the scope of the stability theory, the Hopf bifurcation analysis,
and the investigation of small periodic forcing of stationary solutions of hyperbolic PDEs. The main reason why those
techniques are well established for nonlinear ODEs and parabolic PDEs, but not for nonlinear hyperbolic PDEs, 
is that, in contrast to parabolic case, hyperbolic operators in general do not improve the regularity of their solutions in time
(the question is closely related to propagation of singularities along characteristic curves).
This complicates, in particular, proving the Fredholmness property of the linearizations which is crucial for the analysis of solutions
to nonlinear problems. We provide a range of boundary conditions ensuring the desired smoothingness effect,
which still makes possible to handle the bifurcation analysis of a class of nonlinear problems 
(this idea goes back to \cite{akr}).

The practical motivation is caused by applications to mathematical biology~\cite{hiha}, chemical kinetics (describing mass transition in terms
of convective diffusion and chemical reaction and analysis of chemical processes in counterflow chemical reactors~\cite{AkrBelZel,Sli,Zel1,Zel2}), 
and semiconductor laser dynamics (describing the appearance of
self-pulsations of lasers and modulation of stationary laser states by time periodic electric pumping~\cite{LiRadRe,Rad,RadWu}).

\section{Continuous initial data}\label{sec:classical}
Here we consider the case of continuous initial data $\vphi(x)$.
We will assume the 
 zero-order compatibility conditions
between  (\ref{eq:2}) and (\ref{eq:3}), namely
\begin{equation}\label{eq:4}
\begin{array}{lccl}
\varphi_i(0)&=&h_i(0,v(0)),&\quad k+1\le i\le n,\\
\varphi_i(1)&=&h_i(0,v(0)),&\quad 1\le i\le k,
\end{array}
\end{equation}
where
$
v(0)=\Bigl(\varphi_1(0),\dots,\varphi_k(0),\varphi_{k+1}(1),\dots,
\varphi_{n}(1)\Bigr).
$
By $\|\cdot\|$ we denote the Euclidian norm in $\R^n$.

\begin{theorem}[{\cite[Thm~3.1]{ijdsde}}]\label{thm:ijdsde}
 Assume that the data $\lambda_i$, $a_{ij}$, $g_i$, $\varphi_i$, and $h_i$ are 
continuous functions in all their arguments, and the coefficients $\la_i$ are Lipschitz in $x\in[0,1]$
locally in $t\in[0,\infty)$.
Suppose that $h_i(t,z)$ are continuously differentiable in $z\in\R^n$ and for each $T>0$ 
there exists $C>0$ such that
\begin{equation}\label{eq:h}
\left\|\nabla_zh(t,z)\right\|\le C\left(\log\log H(t,\|z\|)\right)^{1/4},
\end{equation}
where $H$ is a polynomial in $\|z\|$ with
coefficients in $C[0,T]$.
If the zero-order compatibility conditions (\ref{eq:4}) are
fulfilled, then problem (\ref{eq:1})--(\ref{eq:3}) has a unique
continuous solution in $\overline\Pi$ which can be found by the 
sequential approximation method.
\end{theorem}

We now introduce the notions of an {\it Expansion Path} and an 
{\it Influence Path}, that will be our main technical tools. 
Let $\omega_i(\tau;x,t)$ denote
 the  characteristic of the $i$-th equation of (\ref{eq:1})
passing through $(x,t)\in\overline{\Pi}$.
Let $\chi$ be a characteristic of the $i$-th equation of system (\ref{eq:1}).
Suppose that $\chi$ reaches $\d\Pi$ at two points. 
Let $(x,t)$ be that of  these points having larger ordinate (hence, $x=0$ or $x=1$). 
We say that $\chi$ {\it reflects} at $(x,t)$ if 
\begin{equation}\label{eq:ij}
\d_{z}h_j(t,v_1(t),\dots,v_{i-1}(t),z,v_{i+1}(t),\dots,v_n(t))\ne 0 \mbox{ for some } j\le n \mbox{ and } z\in\R.
\end{equation}
In this case that of the characteristics $\om_j(\tau;x,t)$ and $\om_j(\tau;1-x,t)$ which lies above the line $\{\tau=t\}$
is called a {\it reflection} of  $\chi$. If the reflection is defined by $\om_j(\tau;1-x,t)$, it will be called a 
{\it jumping reflection}.

\begin{remark}
Note that condition (\ref{eq:ij}) means that the $i$-th components of the vector $v(t)$
participates in evaluation  of $u_j(0,t)$ for $k+1\le j\le n$ and of $u_j(1,t)$ for $j\le k$.
Whenever the continuous solution
to problem (\ref{eq:1})--(\ref{eq:3}) is known (see Theorem~\ref{thm:ijdsde}),
condition (\ref{eq:ij}) is easily checkable. Otherwise, one can use 
the following constructive sufficient condition: 
(\ref{eq:ij}) holds true whenever for every $(z_1,\dots,z_{i-1},z_{i+1},\dots,z_n)\in\R^{n-1}$
there is $j\le n$ such that
$
\d_{z_i}h_j(t,z_1,\dots,z_n)\not\equiv 0.
$
\end{remark}

\begin{defn}\rm
A sequence of characteristics $\chi_1,\chi_2,\dots$ is called an {\it Extension Path (EP)} if each $\chi_{l+1}$ is a reflection of $\chi_{l}$
(see Figure~1).
\end{defn}

\begin{defn}\label{thm:IP}\rm
A sequence of curved segments $B_1,\dots,B_s$ is called an {\it Influence Path (IP)} if  
the following conditions are met.
\begin{itemize}
\item
Each $B_l$ is a continuous part of a characteristic $\chi_l$ of the $j_l$-th equation for some $j_l\le n$;
\item
The whole path is monotone in the sense that the coordinate $t$ continuously increases while moving along it;
\item
The transition from $B_l$ to $B_{l+1}$ can be of three types:
\begin{enumerate}
\item[--]
$B_l$ and $B_{l+1}$ meet at a point $(x,t)$ such that 
$a_{j_{l+1}j_l}(x,t)\not\equiv 0$ in any  neighborhood of $(x,t)$;
\item[--]
$B_l$ and $B_{l+1}$ meet at a point $(x,t)$  with $x=0$ or $x=1$ and in any neighborhood of $t$ a characteristic of the $j_{l}$-th equation reflects   to a characteristic 
of the $j_{l+1}$-th equation;
\item[--]
$B_l$ terminates at a point $(x,t)$ with $x=0$ or $x=1$, $B_{l+1}$ starts at the point $(1-x,t)$
on the opposite side of $\d\Pi$, and in any neighborhood of $t$ 
a characteristic of the $j_{l}$-th equation 
makes a jumping reflection to a characteristic 
of the $j_{l+1}$-th equation.
\end{enumerate}
\end{itemize}
\end{defn}
Roughly speaking,  an IP is a piecewise continuous curve  with smooth peaces $B_1,\dots,B_s$ lying on 
characteristic curves such that either  $B_l$ and $B_{l+1}$ meet within $\overline\Pi$ or $B_l$ terminates at 
$(x,t)\in\d\Pi$ and $B_{l+1}$ starts at the point $(1-x,t)$ on the opposite side of $\d\Pi$.

\begin{figure}
\centerline{\includegraphics{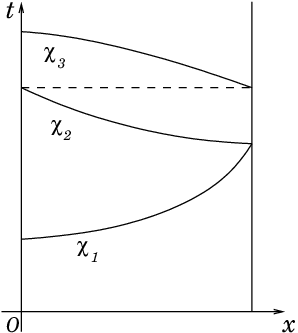}}
\caption{An EP.}
\end{figure}

\begin{remark}
Note that any EP is an IP, while the opposite is not necessary true because segments of an IP
not necessary lie on reflected characteristics.
\end{remark}

\begin{defn}\rm
Define a set $X_i$, called the {\it domain of influence of the initial data on  $u_i$},
as follows: $(x,t)\in X_i$ if  the value  $u_i(x,t)$  can be changed by varying a
 function $\vphi$ in (\ref{eq:2}).
\end{defn}
Since initial data expand inside $\Pi$ along characteristic curves according to  boundary conditions (\ref{eq:3})
and the lower order terms of system (\ref{eq:1}), we have the following characterization.

\begin{lemma}\label{lem:B}
$(x,t)\in X_i$ if and only if there is an IP emanating from the initial axis $t=0$ and going through a part of the characteristic $\om_i(\tau;x,t)$.
\end{lemma}
The sufficiency follows from the proof of the necessity in Theorem~\ref{thm:I}. The proof of the necessity
is based on the constructive description of the domain of dependence of $u_i$ at $(x,t)$ which turns
out to be the union of the IPs going through a part of $ \om_i(\tau;x,t)$ below the line 
$\tau=t$. Under the domain of dependence of $u_i$ at $(x_0,t_0)$ we mean the set of all points $(x^\prime,t^\prime)\in\overline\Pi$
such that if varying the function
\begin{equation}\label{eq:psi}
u(x,t^\prime)=\psi(x)
\end{equation}
in any sufficiently small neighborhood of $x^\prime$, the $i$-th component of the solution to problem 
(\ref{eq:1}), (\ref{eq:psi}), (\ref{eq:3}) changes at $(x_0,t_0)$.

\begin{defn}\label{thm:circ}
\rm
Define a set $X_i^{\circ}$ to be the union of IPs emanating from the initial axis $t=0$ 
and satisfying Definition~\ref{thm:IP} with additional conditions imposed on the three transition 
types from $B_l$ to $B_{l+1}$:
\begin{enumerate}
\item[--]
$B_l$ and $B_{l+1}$ meet at a point $(x,t)$ such that $a_{j_{l+1}j_l}(x,t)\ne 0$;
\item[--]
$B_l$ and $B_{l+1}$ meet at a point $(x,t)$  with $x=0$ or $x=1$ and at this point $\chi_{l}$  reflects to
$\chi_{l+1}$;
\item[--]
$B_l$ terminates at a point $(x,t)$ with $x=0$ or $x=1$, $B_{l+1}$ starts at the point $(1-x,t)$
on the opposite side of $\d\Pi$, and at $(x,t)$ the characteristic $\chi_{l}$ 
makes a jumping reflection to $\chi_{l+1}$.
\end{enumerate}
\end{defn}
Note that $\overline{X_i^{\circ}}=X_i$.

We now introduce two conditions that will occur in formulation of our results.

($\io$) For every  $T>0$ there exists $T'>T$ such that, for all $x\in[0,1]$, every EP
passing through $(x,T)$ lies
below the line $t=T'$.

($\io\io$)
For every $T>0$ and $i\le n$ there exists
$T'>T$ such that, for all $x\in[0,1]$ with $(x,T)\in X_i^{\circ}$, 
every EP containing the characteristic $\om_i(\tau;x,T)$ lies below the line $t=T'$.

\begin{remark}
In many important cases conditions ($\io$) and ($\io\io$) can easily be reformulated and verified in terms of 
$\lambda$, $A$, and $h$. Sometimes they can be verified even directly (see examples in Section~3).
\end{remark}
Set
$
\Pi^T=\{(x,t)\,|\,0<x<1, 0<t<T\}.
$

\begin{remark}
Let $T>0$. In the domain $\Pi\setminus\overline{\Pi^{T}}$
let us consider problem  (\ref{eq:1}), (\ref{eq:psi}), (\ref{eq:3}) with 
$a_{ij}\equiv 0$ for all $i\ne j$
(i.e., system (\ref{eq:1}) is decoupled) and with $t^\prime$ replaced by $T$ in (\ref{eq:psi}). Then
condition $(\io)$ means that, whatever $T>0$ and $\psi(x)$, the function $\psi(x)$ has a bounded domain of
influence on $u_i$ for every $i\le n$. In other words, for any decoupled system (\ref{eq:1}), 
if $\psi(x)$ is singular at some point $x\in[0,1]$, then this singularity 
 expands outside $\overline{\Pi^T}$ along a
finite number of characteristic curves within $\overline{\Pi^{T^\prime}}$ for
some $T^\prime>T$ that does not depend on $x\in[0,1]$.
In contrast to $(\io)$, condition $(\io\io)$ means that, whatever $T>0$ and $\psi_i(x)$, the function $\psi_i(x)$
restricted to $X_i^\circ\cap\{t=T\}$ has a bounded domain of
influence on $u_i$ for every $i\le n$.
\end{remark}

\begin{theorem}\label{thm:I}
 Assume that the data $\lambda_i$, $a_{ij}$, $g_i$, and $h_i$ are 
smooth functions in all their arguments and conditions  (\ref{eq:L1}) and (\ref{eq:h}) are
fulfilled. 
\begin{itemize}
\item
{\rm Sufficiency.} Assume that condition ($\io$) is fulfilled. Then the continuous solution to problem (\ref{eq:1})--(\ref{eq:3}) 
is smoothing for any  $\vphi\in C[0,1]^n$ satisfying equalities (\ref{eq:4}).
\item
{\rm Necessity.} Assume that  the continuous solution to problem (\ref{eq:1})--(\ref{eq:3}) 
is smoothing for any $\vphi\in C[0,1]^n$ satisfying equalities (\ref{eq:4}). Then condition
($\io\io$) is fulfilled.
\end{itemize}
\end{theorem}

\begin{remark}
One can easily see that the sufficient condition ($\io$) implies
the necessary condition ($\io\io$),
while the converse is not true. Nevertheless, in a quite general situation  when 
$X_i^{\circ}=\overline\Pi$ for every $i\le n$, it is not difficult to observe that
conditions $(\io)$ and $(\io\io)$ are equivalent. Hence we have the following result.
\end{remark}

\begin{corollary}
Assume that the data $\lambda_i$, $a_{ij}$, $g_i$, and $h_i$ are 
smooth functions in all their arguments, $X_i^{\circ}=\overline\Pi$ for every $i\le n$,
and conditions  (\ref{eq:L1}) and (\ref{eq:h}) are fulfilled. 
Then the continuous solution to problem (\ref{eq:1})--(\ref{eq:3}) is smoothing for any  $\vphi\in C[0,1]^n$ if and only if 
 condition ($\io$) is fulfilled.
\end{corollary}

\begin{proof}
{\it Sufficiency.}
Assume that condition ($\io$) is fulfilled.
Define a sequence $T_1,T_2,\dots$ inductively by the following rule. Let $T_1$ be the infimum of those $\tau>0$
for which there is $x$ and an EP passing through  $(x,0)$ and lying below the line $t=\tau$; let $T_j$ for $j>1$ be the infimum 
of those $\tau>0$ for which there is $x$ and an EP passing through $(x,T_{j-1})$ and  lying below the line $t=T_{j-1}$.
Note that $T_j$ is monotone and approaches the infinity. 
The latter fact is a simple consequence of the smoothness assumptions on $\Lambda$.
By Theorem~\ref{thm:ijdsde}, problem (\ref{eq:1})--(\ref{eq:3}) has a 
unique continuous solution $u$ in $\overline\Pi$. 
It suffices to show that $u\in C^j\left(\overline\Pi\setminus{\Pi^{T_j}}\right)^n$.

Consider  problem (\ref{eq:1})--(\ref{eq:3}) first in $\overline\Pi\setminus\Pi^{T_1}$. The
solution satisfies the system of integral equations
\begin{eqnarray}
\lefteqn{u_i(x,t) = E_i(t_i(x,t);x,t)h_i(t_i(x,t),v(t_i(x,t)))}\label{eq:u}\\&&
+
\int\limits_{t_i(x,t)}^tE_i(\tau;x,t)\left[
-\sum\limits_{j\ne i}a_{ij}(\xi,\tau)u_{j}(\xi,\tau)+g_{i}(\xi,\tau)
\right]\bigg|_{\xi=\omega_i(\tau;x,t)}\,d\tau,\quad
i\le n,\nonumber
\end{eqnarray}
where
$
E_i(\tau;x,t)=\exp\int_t^\tau a_{ii}(\omega_i(\tau_1;x,t),\tau_1)\,d\tau_1
$
and $t_i(x,t)$  denotes
the smallest value of $\tau\ge 0$ at which the characteristic $\xi=\om_i(\tau;x,t)$
 reaches $\partial\Pi$. In the sequel, along with the equation $\xi=\om_i(\tau;x,t)$
we will also use its inverse form $\tau=\tilde\om_i(\xi;x,t)$. Due to the definition of
$T_1$, in the boundary term $h_i(t_i(x,t),v(t_i(x,t)))$ in (\ref{eq:u}) we
can substitute
\begin{align}
&v_j(t_i(x,t)) = E_j(t_j(x_j,t_i(x,t));x_j,t_i(x,t))h_j(t_j(x_j,t_i(x,t)),v(t_j(x_j,t_i(x,t))))\label{eq:v}\\&
+
\int\limits_{t_j(x_j,t_i(x,t))}^{t_i(x,t)}E_j(\tau;x_j,t_i(x,t))\left[
-\sum\limits_{s\ne j}a_{js}(\xi,\tau)u_{s}(\xi,\tau)+g_{j}(\xi,\tau)
\right]\bigg|_{\xi=\omega_j(\tau;x_j,t_i(x,t))}\,d\tau,
\nonumber
\end{align} 
where $x_j=0$ if $1\le j\le k$ and $x_j=1$ if $k+1\le j\le n$. Continuing in this fashion, 
the right-hand side of (\ref{eq:u}) can  eventually be brought
into a form  depending neither on $v$ nor on
$\vphi$. This version of $(\ref{eq:u})$ will be referred to as $(\ref{eq:u}^\prime)$. We begin with establishing
the $C_{x}^{1}\left(\overline\Pi\setminus\Pi^{T_1}\right)$-smoothness of $u$. It will be proved once we
show that the right-hand side of $(\ref{eq:u}^\prime)$ has a continuous partial derivative in $x$.
The latter can be done by transforming (appropriate changing of variables in) all integrals occurring in $(\ref{eq:u}^\prime)$.
The transformation of each integral follows the same scheme, which we illustrate by example
of the integral expression 
\begin{equation}\label{eq:I0}
I_{ijm}(x,t)=\int\limits_{t_j(x_j,t_i(x,t))}^{t_i(x,t)}E_j(\tau;x_j,t_i(x,t))a_{jm}(\xi,\tau)u_{m}(\xi,\tau)
\bigg|_{\xi=\omega_j(\tau;x_j,t_i(x,t))}\,d\tau.
\end{equation}
We will use assumption (\ref{eq:L1}). Suppose, for instance, that $i\le k$, $j\ge k+1$, $m\le k$, and $t_m(1,t_j(x_i,t_i(x,t)))>0$. 
This entails, in particular, that $x_i=0$ and $x_j=1$. 
\begin{figure}
\centerline{\includegraphics{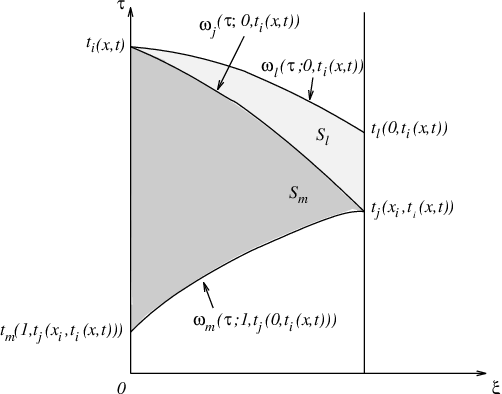}}
\caption{The domain of integration.}
\end{figure}
Due to
(\ref{eq:u}), we obtain (up to the sign)
\begin{align} 
&I_{ijm}(x,t)=\int\limits_{t_j(0,t_i(x,t))}^{t_i(x,t)}E_j(\tau;0,t_i(x,t))a_{jm}(\xi,\tau)
\biggl[
E_m(t_m(\xi,\tau);\xi,\tau)h_m(t_m(\xi,\tau),v(t_m(\xi,\tau)))&
\nonumber\\
&+
\int\limits_{t_m(\xi,\tau)}^{\tau}E_m(\tau_1;\xi,\tau)\biggl\{
-\sum\limits_{p\ne m}(a_{mp}u_{p}
+g_{m})(z,\tau_1)\biggr\}
\bigg|_{z=\omega_m(\tau_1;\xi,\tau)}\,d\tau_1\biggr]\bigg|_{\xi=\omega_j(\tau;0,t_i(x,t))}d\tau&\nonumber
\\&
=\int\limits_{t_m(1,t_j(0,t_i(x,t)))}^{t_i(x,t)}E_j(\rho;0,t_i(x,t))
a_{jm}(\eta,\rho)
E_m(\tau;\eta,\rho)
\frac{\la_m(0,\tau)}{(\la_m-\la_j)(\eta,\rho)}&\nonumber
\\&
\displaystyle
\times\exp\left\{
\int_\tau^\rho\la_{i\xi}^{'}(\xi,\sigma)\bigg|_{\xi=\om_m(\sigma;\eta,\rho)}\,d\sigma
\right\}
\bigg|_{\rho=\theta(0,\tau;x,t),\eta=\omega_m(\rho;0,\tau)}
h_m(\tau,v(\tau))d\tau&
\nonumber
\end{align} 
\begin{align} 
&+\int_{S_m}E_j(\rho;0,t_i(x,t))a_{jm}(\eta,\rho)E_m(\tau;\eta,\rho)
\frac{1}{(\la_m-\la_j)(\eta,\rho)}&\nonumber
\\&
\displaystyle
\times\exp\left\{
\int_\tau^\rho\la_{i\xi}^{'}(\xi,\sigma)\bigg|_{\xi=\om_m(\sigma;\eta,\rho)}\,d\sigma
\right\}
\bigg|_{\rho=\theta(\xi,\tau;x,t),\eta=\omega_m(\rho;\xi,\tau)}&
\nonumber
\\&
\displaystyle
\times\left\{
-\sum\limits_{p\ne m}a_{mp}(\xi,\tau)u_{p}(\xi,\tau)+g_{m}(\xi,\tau)
\right\}
\,d\xi\d\tau,&
\label{eq:I}
\end{align} 
where $S_m$ is the area shown in Fig.~2 and $\theta(\xi,\tau;x,t)$ denotes the $t$-coordinate of the point 
where the characteristics $\om_j(\tau_1;0,t_i(x,t))$ and $\om_m(\tau_1;\xi,\tau)$ intersect. 
The other cases are similar. For example, if $m\ge k+1$, then in the formula (\ref{eq:I})
 index $m$ should be  
replaced by $l$ and the integration over $[t_m(1,t_j(0,t_i(x,t))),t_i(x,t)]$ should be
replaced by the integration over $[t_j(0,t_i(x,t)),t_l(0,t_i(x,t))]$.

The desired $C_x^1$-smoothness of $I_{ijm}(x,t)$ follows from the  $C_x^1$-smoothness of $\theta(\xi,\tau;x,t)$.
The latter is a consequence of the smoothness properties
of $\Lambda$. Indeed, from the equality
$
\om_j(\theta(\xi,\tau;x,t);0,t_i(x,t))=\om_m(\theta(\xi,\tau;x,t);\xi,\tau),
$
we conclude that, if $\d_x\theta$ exists, then it is given by the formula
\begin{align} 
\d_x\theta&(\xi,\tau;x,t)=\frac{\d\om_j(\theta(\xi,\tau;x,t);0,t_i(x,t))}{\d t_i}
\frac{\d t_i(x,t)}{\d x}
\nonumber\\&
\times
\left(\frac{\d\om_m(\theta(\xi,\tau;x,t);\xi,\tau)}{\d \theta}-
\frac{\d\om_j(\theta(\xi,\tau;x,t);0,t_i(x,t))}{\d \theta}\right)^{-1}
\nonumber\\&
=-\la_j(0,t_i(x,t))\e\left\{\int\limits_{t_i(x,t)}^{\theta(\xi,\tau;x,t)}
\la_{j\xi}^{'}(\xi,\sigma)\bigg|_{\xi=\om_j(\sigma;0,t_i(x,t))}\,d\sigma\right\}
\nonumber\\&
\times(\la_i(0,t_i(x,t)))^{-1}\e\left\{\int\limits_{t}^{t_i(x,t)}
\la_{i\xi}^{'}(\xi,\sigma)\bigg|_{\xi=\om_i(\sigma;x,t)}\,d\sigma\right\}
\nonumber\\&
\times\Bigl[\la_m(\om_m(\theta(\xi,\tau;x,t);\xi,\tau),\theta(\xi,\tau;x,t))
\nonumber\\&
-\la_j(\om_j(\theta(\xi,\tau;x,t);0,t_i(x,t)),\theta(\xi,\tau;x,t))\Bigr]^{-1}.
\nonumber
\end{align} 
Thanks to the equality $\om_m(\theta(\xi,\tau;x,t);\xi,\tau)=\om_j(\theta(\xi,\tau;x,t);0,t_i(x,t))$
and condition (\ref{eq:L1}),
the function $I_{ijm}(x,t)$ is continuously differentiable in $x$. 

Thus, the right-hand side of $(\ref{eq:u}^\prime)$ is continuously differentiable in $x$.
Therefore, $\d_xu\in C\left(\overline\Pi\setminus\Pi^{T_1}\right)^n$.
The membership of $u$ in $C^1\left(\overline\Pi\setminus\Pi^{T_1}\right)^n$ now directly follows
from system~(\ref{eq:1}).

In the next step, we prove that $u\in C^2\left(\overline\Pi\setminus\Pi^{T_2}\right)^n$.
For $\d_x u$ we have equations
\begin{equation}\label{eq:dxu}
\begin{array}{ll}
\displaystyle
\partial_xu_i(x,t) = (R_{i}u)(x,t)
+\int\limits_{t_i(x,t)}^tE_i(\tau;x,t)\biggl[
-\sum\limits_{j\ne i}a_{ij}(\xi,\tau)\d_\xi u_{j}(\xi,\tau)
\\
-\partial_\xi\lambda_i(\xi,\tau)\partial_\xi u_i(\xi,\tau)
+\left(\partial_\xi f_i\right)(\xi,\tau,u)\biggr]
\bigg|_{\xi=\omega_i(\tau;x,t)}
\,d\tau,\qquad i\le n,
\end{array}
\end{equation}
where   
$
f_i(\xi,\tau,u)=g_i(\xi,\tau)
-\sum_{j=1}^na_{ij}(\xi,\tau)u_{j},
$
\begin{eqnarray}\label{eq:R}
(R_{i}u)(x,t)&=&E_i(\tau;x,t)
\lambda_i^{-1}(y_i,\tau)\nonumber
\\&&\times
\Bigl[
f_i(y_i,\tau,u)-\nabla_vh_i(\tau,v)\cdot v^\prime(\tau)-
(\partial_\tau h_i)(\tau,v)
\Bigr]\Big|_{\tau=t_i(x,t)},
\end{eqnarray}
$y_i=0$ for $k+1\le i\le n$, and $y_i=1$ for $1\le i\le k$.
Here $w\cdot z$ denotes the scalar product in $\R^n$. 
In (\ref{eq:R}) we can represent $v^\prime(t)$ in the form
\begin{eqnarray}
v_j^\prime(t)&=&f_j(x_j,t,u)-\la_j(x_j,t)\d_xu_j(x_j,t)=
f_j(x_j,t,u)-\la_j(x_j,t)\biggl[ (R_{j}u)(x_j,t)
\nonumber\\&
+&\int\limits_{t_j(x_j,t)}^tE_j(\tau;x_j,t)\biggl(
-\sum\limits_{s\ne j}a_{js}(\xi,\tau)\d_\xi u_{s}(\xi,\tau)
-\partial_\xi\lambda_j(\xi,\tau)\partial_\xi u_j(\xi,\tau)\nonumber\\&
+&\left(\partial_\xi f_j\right)(\xi,\tau,u)\biggr)
\bigg|_{\xi=\omega_j(\tau;x_j,t)}
\,d\tau\biggr].\nonumber
\end{eqnarray}
We continue in this fashion up to getting a representation of the boundary term (\ref{eq:R})
that does not depend on $v^\prime$, what is
possible due to the  definition of $T_2$. To show that the right-hand side of 
the obtained expression for $\partial_xu_i(x,t)$, say $(\ref{eq:dxu}^\prime)$,
is continuously differentiable in $x$, we transform all integrals contributing into
$(\ref{eq:dxu}^\prime)$ similarly to (\ref{eq:I}). Using the fact that
$u\in C^1(\overline\Pi\setminus\Pi^{T_1})^n$, one can easily show that
 $u\in C_{x,t}^{2,1}\left(\overline\Pi\setminus\Pi^{T_2}\right)^n$.
The desired $C^{2}\left(\overline\Pi\setminus\Pi^{T_2}\right)$-smoothness of the solution 
is then a direct consequence of system (\ref{eq:1}) and its differentiations.

We further proceed by induction on $k$.
Assuming that $u\in C^{k-1}\left(\overline\Pi\setminus\Pi^{T_{k-1}}\right)^n$ for some
$k\ge 2$, we will prove that
$u\in C^{k}\left(\overline\Pi\setminus\Pi^{T_{k}}\right)^n$. We differentiate (\ref{eq:1})
$k$ times in $x$,  thereby obtaining a  system for $\d_x^{k}u$. 
The integral form of this system is similar to $(\ref{eq:dxu})$. 
Analogously to $(\ref{eq:dxu})$, the definition of $T_k$ makes possible an integral
representation of the system for $\d_x^{k}u$ which does not depend on $v^{(k)}$
and includes integrals of $\d_x^{k-1}u$ similar to $I_{ijs}(x,t)$.
To show the $C_x^k$-smoothness of the solution, we transform the integral terms analogously to
(\ref{eq:I}). Finally, the $C^k$-smoothness of $u$ outside of $\Pi^{T_k}$
follows from suitable differentiations of system (\ref{eq:1}).
The sufficiency is thereby  proved.

{\it Necessity.}
Suppose that  condition
($\io\io$) is not fulfilled and prove that the solution to problem (\ref{eq:1})--(\ref{eq:3}) 
is not smoothing for some  $\vphi\in C[0,1]^n$.
Fix $T=t_0>0$ and $i\le n$ such that for all $T^\prime>t_0$ there is an EP containing the characteristic
$\om_i(\tau;x,t_0)$ for some  $(x,t_0)\in X_i$ and going beyond  $\overline{\Pi^{T^\prime}}$.

Since all  singularities of solutions expand along EPs, it is sufficient 
to prove that there exist $m\in\N$ and $\vphi\in C[0,1]^n$ such that, whatever $(x_0,t_0)\in X_i^\circ$, the solution $u$ is not 
$C^m$-smooth at $(x_0,t_0)$. 

By Lemma~\ref{lem:B}, for any $(x,t_0)\in X_i$ there is an  Influence Path 
$B_1,\dots,B_s$ emanating from the line $t=0$ and going through a part of $\om_i(t;x,t_0)$.
Denote the smallest possible ``length'' of such a path by $s(x)$. By the smoothness
assumptions on the initial data, if $x^\prime$ is sufficiently close to $x$ and $(x^\prime,t_0)\in X_i$, then $s(x^\prime)\le s(x)$.
Since $X_i\cap\{t=t_0\}$ is closed, the standard compactness argument implies that $s(x)$ is bounded by a constant $s_0$ uniformly over all  $x$.

Let $\varphi$ be a continuous nowhere differentiable function.
Consider an arbitrary $(x_0,t_0)\in X_i^\circ$. Let us fix a shortest Influence Path 
$B_1,\dots,B_s$ from some $(z_0,0)$  to  $(x_0,t_0)$ with 
 $B_s$ lying on $\om_i(\tau;x_0,t_0)$ and with  the  transition from $B_l$ to $B_{l+1}$ as in Definition~\ref{thm:circ}.  
Using the smoothness assumptions on $A$ and $h$, we can suppose that $0<z_0<1$.
We now intend to prove that the solution $u$ is not $C^{s}$-smooth  at $(x_0,t_0)$. Since $s\le s_0$,
this will give us the necessary part of the theorem.
 Let $(z_l,T_l)$ denote the starting  point of $B_{l+1}$.  Suppose that
 $B_l$ is a part of a characteristic $\chi_l$ of the $j_l$-th equation. As 
$\vphi_{j_1}$ is not continuously differentiable at $z_0$,  the function
$u_{j_1}$ is not $C^1$-smooth along $\chi_1$  due to
the definition of a characteristic.
If $\chi_2$ is a reflection of $\chi_1$, then $u_{j_2}$ is not differentiable at
$(z_1,T_1)$ and hence 
 $u_{j_2}$ is not differentiable along  $\chi_2$.
Otherwise, by Definition~\ref{thm:IP}, we have $a_{j_2j_1}(z_1,T_1)\ne 0$ 
and $u_{j_2}$ is not $C^2$-smooth along $\chi_2$,
because the integral of $a_{j_2j_1}u_{j_1}$ along $\chi_2$ in the integral form of the $j_2$-th 
differential equation is not a $C^2$-function and this nonsmoothness cannot be compensated 
by any other summands in the integral representation of our problem.
 It follows that $u_{j_2}$ is not $C^2$-smooth,
in particular, at $(z_2,T_2)$. Continuing in this way, we  
arrive at the conclusion that $u_i$ is not $C^s$-smooth along $\om_i(t;x_0,t_0)$ and hence $u$ is not $C^s$-smooth at $(x_0,t_0)$. 

Going into the details of the above argument,
let us represent $u_i(x_0,t_0)$ in an integral form with an integration over a neighborhood of the
 Influence Path $B_1,\dots,B_s$. We focus on the case of $s=3$ where characteristic  $\chi_3$ is not a reflection of $\chi_2$ and 
 $\chi_2$ is not a reflection of $\chi_1$  (see Fig.~3). A similar or even simpler  argument works as well
for other possible cases.
Extend each $B_l$ to 
$
B_l^\prime=\chi_{l}\cap\overline{\Pi^{T_{l}}}.
$
Given $\eps>0$, let $\mathcal{B}$ denote the union of $B_s^\prime$  and  one-sided neighborhoods of $B_1^\prime,\dots,B_{s-1}^\prime$, 
 bounded by the characteristics $\om_{j_2}(\tau;\om_i(T_2-\eps;x_0,t_0),T_2-\eps)$ and 
$\om_{j_1}(\tau;\om_{j_2}(T_1+\eps;\om_i(T_2-\eps;x_0,t_0),T_2-\eps),T_1+\eps)$. Fix an $\eps>0$ so that 
$\mathcal{B}$ contains neither $(0,0)$ nor $(1,0)$. 

Now we write an integral representation of $u_i(x_0,t_0)$ in terms of $u$ over $\mathcal{B}$.
  We start from the formula (\ref{eq:u})
for $u_i(x_0,t_0)$ and rewrite a part of the integral in the right-hand side as

\begin{figure}
\centerline{\includegraphics{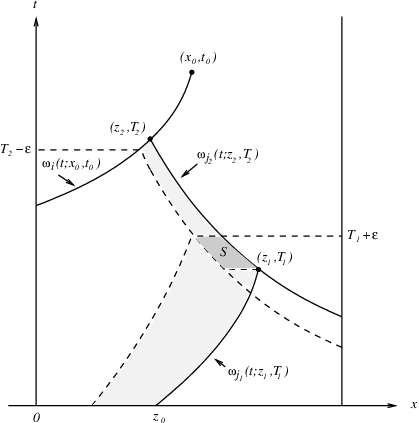}}
\caption{The set $\tilde{\mathcal{B}}$.} 
\end{figure}

\vskip-7mm

\begin{eqnarray}
\lefteqn{
\int\limits_{T_2-\eps}^{T_2}E_i(\tau;x_0,t_0)a_{ij_2}(\xi(\tau),\tau)
u_{j_2}(\xi(\tau),\tau)\,d\tau
}\nonumber\\&&
=\int\limits_{T_2-\eps}^{T_2}E_i(\tau;x_0,t_0)a_{ij_2}(\xi(\tau),\tau)
\biggl[
E_{j_2}(t_{j_2}(\xi(\tau),\tau);\xi(\tau),\tau)
h_{j_2}\left(t_{j_2}(\xi(\tau),\tau),v(t_{j_2}(\xi(\tau),\tau))\right)
\nonumber\\&&
+\int\limits_{t_{j_2}(\xi(\tau),\tau)}^\tau E_{j_2}(\tau_1;\xi(\tau),\tau)
\biggl[
\sum\limits_{s\ne j_2}a_{j_2s}(\om_{j_2}(\tau_1;\xi(\tau),\tau),\tau_1)u_{s}(\om_{j_2}(\tau_1;\xi(\tau),\tau),\tau_1)
\nonumber\\&&
+g_{j_2}(\om_{j_2}(\tau_1;\xi(\tau),\tau),\tau_1)
\biggr]\,d\tau_1\biggr]\,d\tau.
\label{eq:S}
\end{eqnarray}
By construction, $a_{ij_2}(z_2,T_2)\ne 0$. Furthermore, we consider the part of the integral in (\ref{eq:S})
over the area denoted in Fig.~3 by $S$ and transform it similarly to above as follows:
\begin{eqnarray}
\lefteqn{
\int\limits_{T_2-\eps}^{T_2}E_i(\tau;x_0,t_0)a_{ij_2}(\xi(\tau),\tau)
\int\limits_{T_1}^{T_1+\eps}E_{j_2}(\tau_1;\xi(\tau),\tau)
}\nonumber\\&&
\times
a_{j_2j_1}(\omega_{j_2}(\tau_1;\xi(\tau),\tau),\tau_1)
u_{j_1}(\omega_{j_2}(\tau_1;\xi(\tau),\tau),\tau_1)\,d\tau_1\,d\tau
\nonumber\\&&
=\int\limits_{T_2-\eps}^{T_2}E_i(\tau;x_0,t_0)a_{ij_2}(\xi(\tau),\tau)
\int\limits_{T_1}^{T_1+\eps}E_{j_2}(\tau_1;\xi(\tau),\tau)
a_{j_2j_1}(\omega_{j_2}(\tau_1;\xi(\tau),\tau),\tau_1)
\nonumber\\&&
\biggl[\vphi_{j_1}(\omega_{j_1}(0;\omega_{j_2}(\tau_1;\xi(\tau),\tau),\tau_1))
+\int\limits_{0}^{\tau_1}E_{j_1}(\tau_2;\omega_{j_2}(\tau_1;\xi(\tau),\tau),\tau_1)
\nonumber\\&&
\times
\biggl[
\sum\limits_{s\ne j_1}a_{j_1s}(\om_{j_1}\left(\tau_2;\om_{j_2}(\tau_1;\xi(\tau),\tau),\tau_1\right),\tau_2)
u_{s}(\om_{j_1}\left(\tau_2;\om_{j_2}(\tau_1;\xi(\tau),\tau),\tau_1\right),\tau_2)
\nonumber\\&&
+g_{j_1}(\om_{j_1}\left(\tau_2;\om_{j_2}(\tau_1;\xi(\tau),\tau),\tau_1\right),\tau_2)
\biggr]\,d\tau_2
\biggr]\,d\tau_1.
\label{eq:S1}
\end{eqnarray}
Again by construction, $a_{j_2j_1}(z_1,T_1)\ne 0$.
Combining (\ref{eq:u}), (\ref{eq:S}), and (\ref{eq:S1}), we see that
\begin{equation}\label{eq:repr}
u_i(x_0,t_0)=Q_\eps[u](x_0,t_0)+\Phi_i^\eps(x_0,t_0),
\end{equation}
where $Q_\eps$ is a certain operator and
\begin{eqnarray}
\lefteqn{
\Phi_i^\eps(x_0,t_0)=\int\limits_{T_2-\eps}^{T_2}E_i(\tau;x_0,t_0)a_{ij_2}(\xi(\tau),\tau)
\int\limits_{T_1}^{T_1+\eps}E_{j_2}(\tau_1;\xi(\tau),\tau)}
\nonumber\\&&
\times
a_{j_2j_1}(\omega_{j_2}(\tau_1;\xi(\tau),\tau),\tau_1)
\vphi_{j_1}(\omega_{j_1}(0;\omega_{j_2}(\tau_1;\xi(\tau),\tau),\tau_1))\,d\tau_1\,d\tau.
\nonumber
\end{eqnarray}
A similar representation for $u_i(x,t)$ holds
in a sufficiently small neighborhood of $(x_0,t_0)$. From the derivation of (\ref{eq:repr})
it follows  that
$u_i$ at $(x_0,t_0)$ cannot be more regular than any of the
summands in (\ref{eq:repr}). Hence it suffices to show that $\Phi_i^\eps(x_0,t_0)$
is not smooth at $(x_0,t_0)$. Let us make simple transformations:
\begin{eqnarray}
\lefteqn{
\Phi_i^\eps(x_0,t_0)=\int\limits_{T_2-\eps}^{T_2}E_i(\tau;x_0,t_0)a_{ij_2}(\xi(\tau),\tau)
}\nonumber\\&&
\times\int\limits_{\eta(T_1+\eps,\tau)}^{\eta(T_1,\tau)
}
E_{j_2}(\beta(x,\tau);\xi(\tau),\tau)a_{j_2j_1}(\al(x,\tau),\beta(x,\tau))
|J(x,\tau)|\vphi_{j_1}(x)\,dx\,d\tau
\nonumber\\&&
=\int\limits_{T_2-\eps}^{T_2}E_i(\tau;x_0,t_0)a_{ij_2}(\xi(\tau),\tau)
\nonumber\\&&
\times
\biggl[
\int\limits_{0}^{\eta(T_1,\tau)
}
E_{j_2}(\beta(x,\tau);\xi(\tau),\tau)a_{j_2j_1}(\al(x,\tau),\beta(x,\tau))
|J(x,\tau)|\vphi_{j_1}(x)\,dx
\nonumber
\\&&
\displaystyle+
\int\limits_{\eta(T_1+\eps,\tau)
}^0
E_{j_2}(\beta(x,\tau);\xi(\tau),\tau)a_{j_2j_1}(\al(x,\tau),\beta(x,\tau))
|J(x,\tau)|\vphi_{j_1}(x)\,dx
\biggr]\,d\tau
\nonumber\\&&
\displaystyle=
\int\limits_{\eta(T_1,T_2-\eps)
}^{\eta(T_1,T_2)
}
E_i(\rho(y);x_0,t_0)a_{ij_2}(\ga(y),\rho(y))
\nonumber\\&&
\displaystyle\times
\int\limits_0^y
E_{j_2}(\beta(x,\rho(y));\ga(y),\rho(y))a_{j_2j_1}(\al(x,\rho(y)),\beta(x,\rho(y)))
|J(x,\rho(y))|\vphi_{j_1}(x)\,dx\,dy
\nonumber\\&&
\displaystyle+
\int\limits_{\eta(T_1+\eps,T_2-\eps)
}^{\eta(T_1+\eps,T_2)
}
E_i(\zeta(y);x_0,t_0)a_{ij_2}(\sigma(y),\zeta(y))
\nonumber\\&&
\displaystyle\times
\int\limits_y^0E_{j_2}(\beta(x,\zeta(y));\sigma(y),\zeta(y))a_{j_2j_1}(\al(x,\zeta(y)),\beta(x,\zeta(y)))
|J(x,\zeta(y))|\vphi_{j_1}(x)\,dx\,dy,
\nonumber
\end{eqnarray}
where 
$$
J(x,\tau)=\frac{\d(\tau_1,\tau)}{\d(x,\tau)}
$$
denotes the Jacobian of the transformation 
$$
(\tau_1,\tau)\longrightarrow(x,\tau)=
\left(\omega_{j_1}(0;\omega_{j_2}(\tau_1;\xi(\tau),\tau_1),\tau\right),
$$
$\eta(t,\tau)=\omega_{j_1}(0;\omega_{j_2}(t;\xi(\tau),\tau),t)$; $(\al(x,\tau),\beta(x,\tau))$ 
denotes the intersection point of the characteristics $\om_{j_1}(\tau_1;x,0)$ and $\om_{j_2}(\tau_1;\xi(\tau),\tau)$;
$(\ga(y),\rho(y))$ denotes the intersection point of the characteristics $\om_i(\tau;x_0,t_0)$ and $\om_{j_2}(\tau;\om_{j_1}(T_1;y,0),T_1)$;
$(\sigma(y),\zeta(y))$ denotes the intersection point of the $i$-th and the $j_2$-th
characteristics passing through the points $\om_i(\tau;x_0,t_0)$ and $\om_{j_2}(\tau;\om_{j_1}(T_1+\eps;y,0),T_1+\eps)$.
Note that $\al$, $\be$, $\ga$, $\rho$, $\sigma$, and
$\zeta$ are (at least) continuous functions, what easily follows from the smoothness assumptions imposed on the initial data.
It follows that $u_i$ is not $C^3$-smooth at $(x_0,t_0)$ 
(even if $\al$, $\be$, $\ga$, $\rho$, $\sigma$, and $\zeta$ are $C^\infty$-functions).

The necessity is proved.
\end{proof}



\section{Examples}\label{sec:ex}
Here we give some examples to show how the criterion given by Theorem~\ref{thm:I} works. 
Each of these examples is rather general and interesting by its own.
Throughout this section it is supposed that all characteristics of system (\ref{eq:1})
are bounded. This assumption is not restrictive from the practical point of view.
It is true, for example, whenever $\inf|\lambda_i|>0$ for all $i\le n$.

\subsection{Classical boundary conditions}

As a partial case of problem (\ref{eq:1})--(\ref{eq:3}), consider 
(\ref{eq:1}), (\ref{eq:2}) with classical boundary conditions
\begin{eqnarray}
\begin{array}{ll}
u_i(0,t) = h_i(t), & k+1\le i\le n, \ \  t\in(0,\infty)\\
u_i(1,t) = h_i(t), & 1\le i\le k, \quad\quad\ \, t\in(0,\infty).
\end{array}&&\label{eq:3''}
\end{eqnarray}

\begin{theorem}\label{thm:classical}
 Assume that the data $\lambda_i$, $a_{ij}$, $g_i$, and $h_i$ are 
smooth in all their arguments.
Suppose that  condition (\ref{eq:L1}) is
fulfilled. Then the continuous solution to problem (\ref{eq:1}), (\ref{eq:2}), (\ref{eq:3''}) 
is smoothing for any $\vphi\in C[0,1]^n$ satisfying (\ref{eq:4}).
\end{theorem}
This result is a straightforward consequence of Theorem~\ref{thm:I}. Indeed, since no reflection
from the boundary is possible, every EP passing through $(x,t)$ consists of a single characteristic. 
Moreover, due to the smoothness assumptions on $\Lambda$ and the definition of a characteristic, 
for any $T>0$ there is $T^\prime>T$ such that all characteristics passing through the line  $t=T$
lei below $t=T^\prime$.
Hence the sufficient condition $(\io)$
is fulfilled and the theorem follows.

\begin{remark}
Even in the case of smooth classical
 boundary conditions (\ref{eq:3''}),
the domain of influence of the initial data on $u_i$ for every $i\le n$ in general is unbounded
 (due to the lower order terms in (\ref{eq:1})).
In spite of this, the influence of the initial data on the regularity of $u$ becomes
weaker and weaker in time causing the smoothingness effect.
\end{remark}

\subsection{Periodic boundary conditions}
Suppose that at least one component of the solution, say the first one, satisfies periodic boundary condition. 
Specifically, the first equation in (\ref{eq:3}) is written in the form
$
u_1(0,t)=u_1(1,t).
$
Then the domain of influence of initial data (\ref{eq:2}) on $u_1$ is the whole $\overline\Pi$.
Moreover, for each $(x,t)\in\overline\Pi$, there is an unbounded EP passing  through $(x,t)$ (at least one
constructed by means of characteristics of the first equation). This entails that the necessary condition 
$(\io\io)$ is not fulfilled and hence the solution to problem (\ref{eq:1})--(\ref{eq:3})
is not smoothing.

\subsection{Nonseparable linear boundary conditions of local type}

Consider now the reflection boundary conditions
\begin{eqnarray}
\begin{array}{ll}
\displaystyle
u_i(0,t) = \sum\limits_{j=1}^kb_{ij}u_j(0,t), & k+1\le i\le n, \ \  t\in(0,\infty)\\
\displaystyle
u_i(1,t) = \sum\limits_{j=k+1}^nc_{ij}u_j(1,t), & 1\le i\le k, \quad\quad\ \, t\in(0,\infty),
\end{array}&&\label{eq:ck}
\end{eqnarray}
where $b_{ij}$ and $c_{ij}$ are constants. Systems like (\ref{eq:1}), (\ref{eq:2}), (\ref{eq:ck})
cover linearizations of many mathematical models for chemical kinetics~\cite{Zel1,Zel2}. 

Our aim is to reformulate sufficient condition $(\io)$ constructively in terms of data
(specifically, in terms of boundary data) of our problem. Note first that,
if $b_{ij}\ne 0$ for some $k+1\le i\le n$ and $j\le k$ (resp., $c_{ij}\ne 0$
for some $i\le k$ and $k+1\le j\le n$), then  each characteristic of the  $j$-th equation
reflects from the boundary $x=0$ (resp., from the boundary $x=1$), the reflections being characteristics of the $i$-th equation.
Since all EPs passing through $(x,t)$  are constructed by means of subsequent reflections, for each 
EP  passing through $(x,t)$ and consisting of more than one smooth piece, there is a unique sequence
(finite if the EP  passing through $(x,t)$ is bounded and infinite otherwise) either of kind
\begin{equation}\label{eq:bc}
b_{i_2i_1},c_{i_3i_2},b_{i_4i_3},c_{i_5i_4},\dots
\end{equation}
or of kind
\begin{equation}\label{eq:cb}
c_{i_2i_1},b_{i_3i_2},c_{i_4i_3},b_{i_5i_4},\dots
\end{equation}
with nonzero elements.
Condition $(\io)$ is satisfied iff for every $T>0$ all  EPs  passing through $(x,T)$ are bounded uniformly in $x\in[0,1]$. The latter 
is equivalent to the statement that all the sequences (\ref{eq:bc})
and (\ref{eq:cb}) are finite. This means that condition $(\io)$ is
expressible here in the algebraic form, namely as a finite number of equalities of kind
\begin{equation}\label{eq:bcn}
b_{i_2i_1}\cdot c_{i_3i_2}\cdot b_{i_4i_3}\cdot c_{i_5i_4}\cdot\dots\cdot
b_{i_{n-1}i_{n-2}}\cdot c_{i_ni_{n-1}} = 0
\end{equation}
for all possible $i_1,\dots,i_n$ such that the entries of matrices
$B_{n-k,k}=\{b_{ij}\}_{k+1\le i\le n,j\le n}$ and $C_{k,n-k}=\{c_{ij}\}_{i\le k,k+1\le j\le n}$
appear in (\ref{eq:bcn}) no more than once.

Condition (\ref{eq:bcn}) can be easily reformulated in the following form (see also~\cite{LavLyu}). Set
\begin{eqnarray*}
  R&=&\left(
  \begin{array}{cccccc}
    I_{n-k} & B_{n-k,k}\\
    C_{k,n-k} & I_{k}
  \end{array}
  \right),
\end{eqnarray*}
where $I_j$ is a unit matrix of size $j$. Consider an expansion of the determinant of $R$ 
along the first $n-k$ rows, namely, 
\begin{equation}\label{eq:detR}
\det R = \sum\limits_{i_1<\dots<i_{n-k}}D_{i_1,\dots,i_{n-k}}F_{i_1,\dots,i_{n-k}}.
\end{equation}
Here the respective determinants of order $n-k$ are denoted by $D_{i_1,\dots,i_{n-k}}$
and their cofactors  by $F_{i_1,\dots,i_{n-k}}$. 
It turns out that each summand in (\ref{eq:detR}) starting from the second one is a product of kind 
(\ref{eq:bcn}). One can easily show that (\ref{eq:bcn}) is fulfilled iff all summands in 
(\ref{eq:detR}) excepting  for the first one, which is obviously equals to 1, vanish. Note that the latter condition is obtained in~\cite{LavLyu}
for problem   (\ref{eq:1}), (\ref{eq:2}), (\ref{eq:ck}) with coefficients in  (\ref{eq:1}) depending only on $x$.

\subsection{Smoothingness phenomenon for IBVPs for the wave equation}

Using Theorems~\ref{thm:ijdsde} and~\ref{thm:I}, one can immediately obtain correctly
posed IBVPs for the (nonhomogeneous) wave equation
\begin{equation}\label{eq:w1}
(\d_t^2-a^2\d_x^2)u=f(x,t)
\end{equation}
and state a smoothingness result for them. For instance, let us consider (\ref{eq:w1})
subjected to initial conditions
\begin{equation}\label{eq:w2}
u(x,0) = \varphi(x),\quad \d_tu(x,0) = \psi(x)
\end{equation}
and boundary conditions
\begin{eqnarray}
\begin{array}{lcl}
\displaystyle
u(0,t)& =& h_1\left(t,u(1,t),(\d_tu+a\d_xu)|_{x=0}\right), \\
\displaystyle
(\d_tu+a\d_xu)|_{x=1}& =& h_2\left(t,u(1,t),(\d_tu+a\d_xu)|_{x=0}\right).
\end{array}&&\label{eq:w3}
\end{eqnarray}
The problem (\ref{eq:w1})--(\ref{eq:w3}) is equivalent to the following problem for the
$(2\times 2)$-first-order hyperbolic system
\begin{equation}\label{eq:w4}
\displaystyle
(\d_t+a\d_x)u = w, \quad
(\d_t-a\d_x)w = f(x,t)
\end{equation}
\begin{align}
\displaystyle
u(x,0) = \vphi(x), \quad
w(x,0) = \psi(x)+a\frac{d\vphi(x)}{dx}
\label{eq:w5}\\
\begin{array}{ll}
\displaystyle
u(0,t) = h_1\left(t,u(1,t),w(0,t)\right), \\
\displaystyle
w(1,t) = h_2\left(t,u(1,t),w(0,t)\right).
\end{array}&&\label{eq:w6}
\end{align}
Under a continuous solution to problem (\ref{eq:w1})--(\ref{eq:w3}) we mean the first component of the
continuous solution $(u,w)$ to problem (\ref{eq:w4})--(\ref{eq:w6}).  
In the next theorem conditions $(\io)$  and  $(\io\io)$ are supposed to be  formulated correspondingly to problem
(\ref{eq:w4})--(\ref{eq:w6}). The following smoothingness result for
(\ref{eq:w1})--(\ref{eq:w3}) is a straightforward consequence of Theorems~\ref{thm:ijdsde} and~\ref{thm:I}.

\begin{theorem}\label{thm:Iwave}
 Assume that the data $f$ and $h_i$ are 
smooth functions in all their arguments and condition (\ref{eq:h}) is
fulfilled. 
\begin{itemize}
\item
{\rm Sufficiency.} If 
condition  ($\io$) is true, then the continuous solution to problem (\ref{eq:w1})--(\ref{eq:w3}) 
is smoothing for any  $\vphi\in C^1[0,1]$ and  $\psi\in C[0,1]$  satisfying 
the zero-order compatibility conditions between (\ref{eq:w5}) and (\ref{eq:w6}).
\item
{\rm Necessity.} Assume that  the continuous solution to problem (\ref{eq:w1})--(\ref{eq:w3}) 
is smoothing for any  $\vphi\in C^1[0,1]$ and  $\psi\in C[0,1]$ satisfying
the zero-order compatibility conditions between (\ref{eq:w5}) and (\ref{eq:w6}).
Then condition ($\io\io$) is fulfilled.
\end{itemize}
\end{theorem}

\section{Distributional initial data}\label{sec:nonclassical}
Now we address the case of
 distributional initial data. 
Specifically, we investigate system (\ref{eq:1}) with initial conditions
\begin{eqnarray}
\begin{array}{lc}
\displaystyle
u_i(x,0) = b_i(x)\equiv b_{ir}(x)+b_{is}(x)=b_{ir}(x)+\sum\limits_{j=1}^{m_i}c_{ij}\delta^{(l_{ij})}(x-x_{ij}^*), \\
\ \  i\le n, \ \  0<x_{i1}^*<\dots<x_{im_i}^*<1,\quad x\in(0,1)
\end{array}\label{eq:2'}
\end{eqnarray}
and boundary conditions
\begin{eqnarray}
\begin{array}{ll}
\displaystyle
u_i(0,t) = \sum\limits_{j=1}^{n}p_{ij}(t)v_j^u(t)+r_i(t,v^u(t)), & k+1\le i\le n, \ \  t\in(0,\infty)\\
\displaystyle
u_i(1,t) = \sum\limits_{j=1}^{n}p_{ij}(t)v_j^u(t)+r_i(t,v^u(t)), & 1\le i\le k, \quad\quad\ \, t\in(0,\infty),
\end{array}\label{eq:3'}
\end{eqnarray}
where the regular part $b_r(x)=(b_{1r},\dots,b_{nr})$ of the initial data $b(x)=(b_{1},\dots,b_{n})$ 
is a continuous vector-function on $[0,1]$,
$r_i(t,y)\in C^\infty\left([0,\infty)\times\R^n\right)$, $r_i$ and $\nabla_yr_i$ are
bounded in $y\in\R^n$ uniformly in $t$ over any compact subset of $\R_+$, and $p_{ij}(t)\in C^\infty[0,\infty)$.
Also, $c_{ij}\in\R$, $m_i\in\N$, and $l_{ij}\in\N_0$. 
Furthermore,  $v^u(t)$ just redenotes the vector-function $v(t)$ introduced in Section~1 by~(\ref{eq:vv}). 
We will use such notation below with other superscripts.
Without restriction of generality we can suppose that the zero-order compatibility 
conditions between (\ref{eq:2'}) and (\ref{eq:3'})
are fulfilled, namely,
\begin{equation}\label{eq:zero}
\begin{array}{cc}
 \displaystyle
b_{ri}(0)=\sum\limits_{j=1}^{n}p_{ij}(0)v_j^u(0)+r_i(0,v^u(0)),\quad  k+1\le i\le n\\
 \displaystyle
b_{ri}(1)=\sum\limits_{j=1}^{n}p_{ij}(0)v_j^u(0)+r_i(0,v^u(0)),\quad  1\le i\le k,
\end{array}
\end{equation}
where $v^u(0)=\left(b_{r1}(0),\dots,b_{rk}(0),b_{r,k+1}(1),\dots,b_{rn}(1)\right)$.

Note that  system (\ref{eq:1}), (\ref{eq:2'}), (\ref{eq:3'}) 
is unsolvable in the classical sense, since nonlinearities in the boundary conditions and
singularities in the initial data entail nonlinear operations in the space of distributions.
We, therefore, aim at constructing a delta wave solution to (\ref{eq:1}), (\ref{eq:2'}), (\ref{eq:3'}).
For this purpose we regularize the singular part $b_s(x)$ of initial data (\ref{eq:2'})
by means of the convolution with the so-called model delta nets $\left(\vphi_\eps\right)_{\eps>0}$
where $\vphi\in\D(\R)$, $\int\vphi\,dx=1$,  $\vphi_\eps(x)=1/\eps\vphi(x/\eps)$,
 and $\int|\vphi_\eps|\,dx$ is bounded by a constant independent of $\eps$. 
The regularized problem for $u^\eps$ is this:
\begin{equation}\label{eq:ueps}
\begin{array}{rcll}\displaystyle
(\partial_t  + \Lambda\partial_x + A)
 u^\eps &=& g(x,t)\\
\displaystyle
u_i^\eps(x,0)&=&\displaystyle b_{ir}(x)+b_{is}^\eps(x) = b_{ir}(x)+\sum\limits_{j=1}^{m_i}c_{ij}\vphi_\eps^{(l_{ij})}(x-x_{ij}^*),
& 1\le i\le n \\
\displaystyle
u_i^\eps(0,t)&=& \displaystyle\sum\limits_{j=1}^{n}p_{ij}(t)v_j^{u,\eps}(t)+r_i(t,v^{u,\eps}(t)), & k+1\le i\le n\\
\displaystyle
u_i^\eps(1,t)&=& \displaystyle\sum\limits_{j=1}^{n}p_{ij}(t)v_j^{u,\eps}(t)+r_i(t,v^{u,\eps}(t)), & 1\le i\le k.
\end{array}
\end{equation}
Here
$
v^{u,\eps}(t)=(u_1^\eps(0,t),\dots,u_{k}^\eps(0,t),u_{k+1}^\eps(1,t),\dots,u_{n}^\eps(1,t)).
$
Decompose $A=D+F$ into the diagonal and off-diagonal parts $D$ and $F$, respectively. Our aim is now to show
that $u^\eps\to_{\eps\to 0}z+w$ in $\D'(\Pi)$ where a singular part  $z$ and a regular part $w$ are solutions to systems,
respectively, 
\begin{equation}\label{eq:z}
\begin{array}{rcll}
(\partial_t  + \Lambda\partial_x + D)
 z &=& 0\\
z(x,0)&=& b_{s}(x)\\
\displaystyle
z_i(0,t)&=& \displaystyle \sum\limits_{j=1}^{n}p_{ij}(t)v_{j}^z(t), & k+1\le i\le n\\
\displaystyle
z_i(1,t)&=&  \displaystyle\sum\limits_{j=1}^{n}p_{ij}(t)v_{j}^z(t), & 1\le i\le k
\end{array}
\end{equation}
and
\begin{equation}\label{eq:w}
\begin{array}{rcll}
(\partial_t  + \Lambda\partial_x + A)
 w + Fz &=& g(x,t)\\
w(x,0)&=& b_{r}(x)\\
\displaystyle
w_i(0,t)&=&  \displaystyle\sum\limits_{j=1}^{n}p_{ij}(t)v_{j}^w(t) + r_i(t,v^w), & k+1\le i\le n\\
\displaystyle
w_i(1,t)&=&  \displaystyle\sum\limits_{j=1}^{n}p_{ij}(t)v_{j}^w(t) + r_i(t,v^w), & 1\le i\le k
\end{array}
\end{equation}
(the two systems are called {\it nonlinear splitting}).
System (\ref{eq:z}) is linear, corresponds to the singular part of the initial problem, and is responsible 
for propagation of singularities. System (\ref{eq:w}) is nonlinear and corresponds to the regular part of our problem.
Note that, since $F$ is the off-diagonal part of $A$, the singular forcing term $Fz$ contributes into the hyperbolic system
(\ref{eq:w}) transversely to the characteristic directions and, therefore, cannot produce strong singularities for $w$.
We now show that problems (\ref{eq:z}) and (\ref{eq:w}) are uniquely solvable.

\begin{defn}\rm
By $J_*$ we denote the union, over $i\le n$  and $j\le m_i$, of all EPs containing the characteristic $\om_i(\tau;x_{ij}^*,0)$.
Furthermore, $J$ will denote the union, over $i\le n$  and $j\le m_i$, of all EPs emanating from the point $(x_{ij}^*,0)$.
\end{defn}
Note that $J_*\subset J$.
\begin{lemma} 
System  (\ref{eq:z}) has a unique distributional solution $z$ representable as the 
sum of singular distributions concentrated on characteristic curves in the set~$J_*$.
\end{lemma}
The lemma can be easily proved by the method of characteristics.

Let $C\left(\overline\Pi\setminus J\right)$ be the space of piecewise continuous functions on $\overline\Pi$ 
with (possibly) first-order discontinuities on $J$. 
\begin{lemma} 
System (\ref{eq:w}) has a unique  $C\left(\overline\Pi\setminus J\right)^n$-solution $w$.
\end{lemma}
\begin{proof}
The uniqueness can be  easily proved
by considering the problem for the difference of two  $C\left(\overline\Pi\setminus J\right)^n$-solutions
to (\ref{eq:w}) and applying the method of characteristics.

To prove existence, consider $w=\bar w + \tilde w$, where $\bar w$ is the solution to the linear system with the strongly singular
forcing term
\begin{equation}\label{eq:barw}
\begin{array}{rcll}
(\partial_t  + \Lambda\partial_x + D)
 \bar w + Fz &=& g(x,t)\\
\bar w(x,0)&=& b_{r}(x)\\
\displaystyle
\bar w_i(0,t)&=&  \displaystyle\sum\limits_{j=1}^{n}p_{ij}(t)v_{j}^{\bar w}(t), & k+1\le i\le n\\
\displaystyle
\bar w_i(1,t)&=&  \displaystyle\sum\limits_{j=1}^{n}p_{ij}(t)v_{j}^{\bar w}(t), & 1\le i\le k
\end{array}
\end{equation}
and $\tilde w$ is the solution to the nonlinear system
\begin{equation}\label{eq:tildew}
\begin{array}{rcll}
(\partial_t  + \Lambda\partial_x + A)
 \tilde w + F\bar w &=& 0\\
\tilde w(x,0)&=& 0\\
\displaystyle
\tilde w_i(0,t)&=&  \displaystyle\sum\limits_{j=1}^{n}p_{ij}(t)v_{j}^{\tilde w}(t) + r_i(t,v^{\bar w}+v^{\tilde w}), & k+1\le i\le n\\
\displaystyle
\tilde w_i(1,t)&=& \displaystyle \sum\limits_{j=1}^{n}p_{ij}(t)v_{j}^{\tilde w}(t) + r_i(t,v^{\bar w}+v^{\tilde w}), & 1\le i\le k.
\end{array}
\end{equation}
Note that, if $\bar w$ is a solution to (\ref{eq:barw}) and $\tilde w$ is a solution to (\ref{eq:tildew}),
then $\bar w+\tilde w$ is a solution to (\ref{eq:w}) indeed.

To show that system (\ref{eq:barw}) has a (unique) solution 
$\bar w\in C\left(\overline\Pi\setminus J\right)^n$, we use 
the method of characteristics again, thereby representing $\bar w$ explicitly
in the integral form. Since the singularities in the forcing term $Fz$ are strong and 
 transverse to the corresponding characteristics of the hyperbolic system,
the Volterra integrals of $Fz$  in the integral representation of $\overline w$ are 
at worst discontinuous functions
with the first order discontinuities (if any) on~$J$. 

To show that problem (\ref{eq:tildew}) has a (unique) solution 
$\tilde w\in C\left(\overline\Pi\setminus J\right)^n$, we follow the proof 
of Theorem~\cite[Thm~2.1]{ijdsde} with minor changes.
Fix an arbitrary $T>0$ and rewrite (\ref{eq:tildew}) in $\Pi^T$ in an equivalent 
integral form  with the boundary term not depending on $v^{\tilde w}$
(applying a finite number of integrations along characteristic curves up to reaching the initial axis). 
Using  the assumption
that $\nabla_yr_i(t,y)$ is bounded locally in $t\in\R_+$ and globally in $y\in\R^n$,
it is not difficult to check that the operator defined by the right-hand side of the integral
form  of (\ref{eq:tildew}) maps $C\left(\overline\Pi\setminus J\right)^n$ 
into itself and is  contractible on $\overline{\Pi^{\theta_0}}$ for some  $\theta_0\le T$.
The local existence-uniqueness
result in $\overline{\Pi^{\theta_0}}$ then follows from the Banach fixed point theorem.
If $\theta_0<T$, we consider problem (\ref{eq:tildew}) on $\Pi^T\setminus\overline{\Pi^{\theta_0}}$
with the initial condition replaced by $u(x,\theta_0)=\psi(x)\in C\left([0,1]\setminus J\right)^n$, 
where the function $\psi(x)$ is determined in the first step.
 Similarly to the above, the right-hand side of the integral form of this problem defines the operator mapping
$C\left(\left(\overline\Pi\setminus\Pi^{\theta_0}\right)\setminus J\right)^n$ into itself.
Set $\Pi^0=\emptyset$.
Since the contraction property of the operator under consideration does not depend on the initial conditions, 
one can choose the value of
$\theta_0$ from the very beginning so small that the operator 
is contractible on $\overline{\Pi^{(s+1)\theta_0}}\setminus\Pi^{s\theta_0}$
for any $0\le s\le\lceil T/\theta_0\rceil-1$. Applying the Banach fixed point theorem with respect to the 
 domain $\overline{\Pi^{2\theta_0}}\setminus\Pi^{\theta_0}$, we conclude that problem (\ref{eq:tildew})
has  a unique $C\left(\left(\overline\Pi^{2\theta_0}\setminus\Pi^{\theta_0}\right)\setminus J\right)^n$-solution.
The proof of the unique solvability
in $\Pi^{(s+1)\theta_0}\setminus\Pi^{s\theta_0}$ for each
$s\in\{2,\dots,\lceil T/\theta_0\rceil-1\}$ goes over the same argument. 
Since $T>0$ is arbitrary, the  
unique solvability of (\ref{eq:tildew}) in $C\left(\overline\Pi\setminus J\right)^n$ follows. 

We have thus proved that (\ref{eq:w}) has a unique solution $w$  in $C\left(\overline\Pi\setminus J\right)^n$,
which equals to $\bar w+\tilde w$.
\end{proof}

We are now  prepared to state our result about delta waves.
\begin{theorem}\label{thm:delta}
 Assume that the data $\lambda_i$, $a_{ij}$, $g_i$, $p_{ij}$, $r_i$, and $b_{ir}$ are 
continuous functions in all their arguments, $\la_i$ are Lipschitz in $x\in[0,1]$ locally
in $t\in[0,\infty)$, and both
$r_i$ and $\nabla_yr_i$ are in $\L^\infty\left((0,T)\times\R^n\right)$ for any $T>0$. 
 Given $\eps>0$, let
$u^\eps$, $z^\eps$, and $w^\eps$ be the continuous solutions to problems, respectively
(\ref{eq:ueps}), (\ref{eq:z}), and (\ref{eq:w}) with
$z^\eps$, $w^\eps$, and $b_s^\eps$ in place of  $z$, $w$, and $b_s$, respectively.
Then
\begin{enumerate}
\item[{\bf 1.}]
$u^\eps-z^\eps-w^\eps\to 0 \mbox{ in } \L_{loc}^1(\Pi) \mbox{ as } \eps\to 0$;
\item[{\bf 2.}]
$z^\eps\to z  \mbox{ in } \E'(\Pi) \mbox{ as } \eps\to 0$;
\item[{\bf 3.}]
$w^\eps\to w$ in $C(K)$  as  $\eps\to 0$ for any compact subset $K\subset\overline\Pi\setminus J$;
\item[{\bf 4.}]
$u^\eps\to z+w  \mbox{ in } \D'(\Pi) \mbox{ as } \eps\to 0$,
\end{enumerate}
where $z\in\E'(\Pi)^n$ and $w\in C\left(\overline\Pi\setminus J\right)^n$  are the solutions to problems (\ref{eq:z}) and (\ref{eq:w}),
respectively.
\end{theorem}

\begin{proof}
{\bf 1.}
Let $J_*^\eps$ be the union, over $i\le n$, $j\le m_i$, and $\al\in[-\eps;\eps]$, of all EPs 
containing  $\om_i(\tau;x_{ij}^*+\al,0)$. Note that the set $J_*^\eps$ consists of the ``tubes'' forming the support of $z^\eps$.
Set 
$
N^\eps = J_*^\eps\cap\left(\d\Pi\setminus\{t=0\}\right).
$

To show that the splitting of $u^\eps$ into 
a ``regular'' and a ``singular'' parts $w^\eps$ and $z^\eps$ is correct (what is stated in Item 1 of the theorem), 
we consider the following problem for $\al^\eps=u^\eps-z^\eps-w^\eps$:
\begin{equation}\label{eq:r1}
\begin{array}{rcll}
(\partial_t  + \Lambda\partial_x + A)
\al^\eps &=& 0\\
\al^\eps(x,0) &=& 0\\
\displaystyle
\al_i^\eps(0,t) &=& \displaystyle \sum\limits_{j=1}^{n}p_{ij}(t)v_j^{\al,\eps}(t)+r_i(t,v^{u,\eps})-r_i(t,v^{w,\eps}), & k+1\le i\le n\\
\displaystyle
\al_i^\eps(1,t) &=&\displaystyle \sum\limits_{j=1}^{n}p_{ij}(t)v_j^{\al,\eps}(t)+r_i(t,v^{u,\eps})-r_i(t,v^{w,\eps}), & 1\le i\le k.
\end{array}
\end{equation}
Rewrite
\begin{eqnarray}
\lefteqn{
r_i(t,v^{u,\eps})-r_i(t,v^{w,\eps})
}
\nonumber\\&&
=\left(r_i(t,v^{\al+w+z,\eps})-r_i(t,v^{w+z,\eps})\right)+\left(r_i(t,v^{w+z,\eps})-r_i(t,v^{w,\eps})\right)
\nonumber\\&&
=\int\limits_0^1\left(\nabla_yr_i\right)(t,\sigma v^{\al+w+z,\eps}+(1-\sigma)v^{w+z,\eps})\,d\sigma\cdot v^{\al,\eps}
+\left(r_i(t,v^{w+z,\eps})-r_i(t,v^{w,\eps})\right)\nonumber
\\&&
\equiv \sum\limits_{j=1}^{n}K_{ij}^{\eps}(t)v_{j}^{\al,\eps}+ L_{i}^{\eps}(t),
\end{eqnarray}
where the symbol $\equiv$ introduces the short-hand notation $K_{ij}^{\eps}$ and $L_{i}^{\eps}$.
By assumption, $r_i(t,y)$ and $\nabla_yr_i(t,y)$ are in $\L^\infty\left((0,T)\times\R^n\right)$ for any $T>0$.
Hence $K_{ij}^{\eps}(t)$ and $L_{i}^{\eps}(t)$ are bounded uniformly in $t\in[0,T]$ and $\eps>0$
for an arbitrarily fixed $T>0$.

To prove the desired statement, we will use dominated convergence. It is sufficient to show that,
given $T>0$, the sequence $\al^\eps$ is bounded on $\overline\Pi^T$ uniformly in 
$\eps>0$ and $\al^\eps(x,t)\to 0$ as $\eps\to 0$ pointwise off $J_*$.

Fix an arbitrary $T>0$ and rewrite (\ref{eq:r1}) in $\overline\Pi^T$ in the integral form
\begin{equation}\label{eq:r2}
\al_i^\eps(x,t) = (R_i\al^\eps)(x,t)
-
\int\limits_{t_i(x,t)}^t\sum\limits_{j=1}^n\left(a_{ij}\al_j^\eps\right)(\omega_i(\tau;x,t),\tau)\,d\tau,
\quad i\le n,
\end{equation}
where $(R_i\al^\eps)(x,t) = 0$ if $t_i(x,t)=0$ and
\begin{equation}\label{eq:r3}
(R_i\al^\eps)(x,t)
=\sum\limits_{j=1}^{n}\left(p_{ij} + K_{ij}^{\eps}\right)\left(t_i(x,t)\right)v_{j}^{\al,\eps}\left(t_i(x,t)\right)
+ L_{i}^{\eps}(t_i(x,t))
\end{equation}
otherwise. We will use a modification of (\ref{eq:r2}) where the boundary term $(R_i\al^\eps)(x,t)$
does not depend on $v^{\al,\eps}$. More specifically, using (\ref{eq:r2}), we express
 $v_{j}^{\al,\eps}\left(t_i(x,t)\right)$ contributing into (\ref{eq:r3}) in an integral form. 
If the resulting expression
for $\al^\eps(x,t)$ still depends on $v^{\al,\eps}$, we apply the same procedure again.
In a finite number of steps we reach the initial axis and thereby obtain the desired
 representation for $\al_j^\eps(x,t)$. 
It is the sum of a continuous function (of $K_{ij}^{\eps}$, $L_{i}^{\eps}$, and $p_{ij}$) 
and a finite number of integrals over characteristic curves in the 
EPs passing through $\om_i(\tau;x,t)$ and restricted to $\overline{\Pi^t}$. 
Let $N$ be the maximal number of such integrals 
where the range of maximization is $i\le n$ and $(x,t)\in\overline{\Pi^T}$.
Applying the sequential approximation method, we easily derive the apriori estimate in $\overline{\Pi^T}$:
\begin{eqnarray*}
\lefteqn{
\sup\limits_{\eps,i,x,t}|\al_i^\eps|\le \sup\limits_{i,t,\eps}|L_{i}^{\eps}|
\sum\limits_{s=0}^Nn^s\left(\sup\limits_{i,j,t}|p_{ij}|+\sup\limits_{\eps,i,j,t}|K_{ij}^{\eps}|\right)^s
}\\&&
\times\exp\left\{nT\sup\limits_{i,j,x,t}|a_{ij}|
\sum\limits_{s=0}^Nn^s\left(\sup\limits_{i,j,t}|p_{ij}|+\sup\limits_{\eps,i,j,t}|K_{ij}^{\eps}|\right)^s\right\}.
\end{eqnarray*}
This implies that
 $\al^\eps(x,t)$ is bounded uniformly in $(x,t)\in\overline\Pi^T$ and $\eps>0$, as desired.

It remains to prove that $\al^\eps(x,t)\to 0$ as $\eps\to 0$ pointwise off $J_*$. Note that, due
to the support properties of $z^\eps$, we have $v^{w+z,\eps}(t)\equiv v^{w,\eps}(t)$ and hence
$L_i^{\eps}(t)\equiv 0$ on $\d\Pi\setminus(N^\eps\cup\{t=0\})$.
As $N^{\eps_1}\subset N^{\eps_2}$ for all $\eps_1<\eps_2$, we have the identity
$L_i^{\eps}(t)\equiv 0$ on $\d\Pi\setminus(N^{\eps_0}\cup\{t=0\})$
for  an arbitrarily fixed  $\eps_0>0$ and all $\eps\le\eps_0$. Fix $(x,t)\in\overline{\Pi^T}\setminus J_*$
and $\eps_0>0$ so that $(x,t)\in\overline{\Pi^T}\setminus J_*^{\eps_0}$. Then 
$\al^\eps(x,t)$ is representable in the form (\ref{eq:r2}) with $L_{i}^{\eps}(t_i(x,t)) = 0$.
Note that $t_i(x,t)\in\d\Pi\setminus N^{\eps_0}$; hence $v_j^{\al,\eps}(t_i(x,t))$
is expressible in the form (\ref{eq:r2}) again with $L_{j}^{\eps} = 0$.
 Continuing in this way, similarly to the above we arrive at the integral
representation of $\al^\eps$ where the boundary term depends neither on $L^{\eps}$ nor on
$v^{\al,\eps}$. Roughly speaking, each $\al_i^\eps(x,t)$ is a finite sum of integrals 
(of kind as in (\ref{eq:r2})) over characteristics in the 
EPs passing through $\om_i(\tau;x,t)$ and restricted to $\overline{\Pi^t}$.
We split each of the integrals into two parts, the splitting procedure being
illustrated by example of the integral in (\ref{eq:r2}):
\begin{equation}
\int\limits_{t_i(x,t)}^t\,d\tau = \int\limits_{[t_i(x,t),t]\setminus\triangle_{\eps_0}(x,t)}\,d\tau
+ \int\limits_{\triangle_{\eps_0}(x,t)}\,d\tau,\nonumber
\end{equation}
where $\triangle_{\eps_0}(x,t)=\{\tau\in[t_i(x,t),t]\,:\,(\om_i(\tau;x,t),\tau)\in J_*^{\eps_0}\}$. 
Since $\al^\eps$ is  bounded on  $\overline\Pi^T$ uniformly in $\eps>0$, there is $C>0$
(depending on $T$) such that the absolute value of the second integral is bounded from 
above by $C\eps_0$. This gives us the following apriori estimate  in $\overline{\Pi^T}$:
\begin{equation}\label{eq:r4} 
\begin{array}{ccc}
\displaystyle
\sup\limits_{\eps\le\eps_0,i,x,t}|\al_i^\eps|\le C\eps_0
\sum\limits_{s=0}^Nn^s\left(\sup\limits_{i,j,t}|p_{ij}|+\sup\limits_{\eps>0i,j,t}|K_{ij}^{\eps}|\right)^s
\\
\times\displaystyle
\exp\left\{nT\sup\limits_{i,j,x,t}|a_{ij}|
\sum\limits_{s=0}^Nn^s\left(\sup\limits_{i,j,t}|p_{ij}|+\sup\limits_{\eps>0,i,j,t}|K_{ij}^{\eps}|\right)^s\right\}.
\end{array}
\end{equation}
Note that (\ref{eq:r4}) is true with $\eps_1$ in place of $\eps_0$ for an
arbitrarily fixed $\eps_1\le\eps_0$. The desired convergence is thereby proved.
The first statement of the theorem now follows by dominated convergence.

{\bf 2.}
This part  easily follows from diagonality of matrix $D$ by the
method of characteristics.

{\bf 3.}
Given $\eps$, denote by $\bar w^\eps$ the continuous solution to (\ref{eq:barw}) with $z^\eps$ in place of $z$ 
and let $\bar w$ denote the $C\left(\overline\Pi\setminus J\right)^n$-solution to (\ref{eq:barw}).
Let $K$ be an arbitrary compact  and connected subset of $\overline\Pi\setminus J$.
\begin{claim}
$\bar w^\eps\to \bar w$  in $C(K)$ as $\eps\to 0$.
\end{claim}
\begin{subproof}
The difference $\bar w-\bar w^\eps$ satisfies system
\begin{equation}\label{eq:bar-w}
\begin{array}{rcll}
(\partial_t  + \Lambda\partial_x + D)
\left(\bar w^\eps-\bar w\right) + F(z^\eps-z) &=& 0\\
\bar w^\eps(x,0)-\bar w(x,0)&=& 0\\
\displaystyle
\bar w_i^\eps(0,t)-\bar w_i(0,t) &=&  \displaystyle\sum\limits_{j=1}^{n}p_{ij}(t)\left(v_{j}^{\bar w,\eps}(t)-v_{j}^{\bar w}(t)\right), 
& k+1\le i\le n\\
\displaystyle
\bar w_i^\eps(1,t)-\bar w_i(1,t) &=& \displaystyle \sum\limits_{j=1}^{n}p_{ij}(t)\left(v_{j}^{\bar w,\eps}(t)-v_{j}^{\bar w}(t)\right), 
& 1\le i\le k.
\end{array}
\end{equation}
Since $z^\eps\to z$ in $\E'(\Pi)$ as $\eps\to 0$, we obtain the convergence $\bar w^\eps\to \bar w$ in $\E'(\Pi)$ as $\eps\to 0$.
It remains to prove that $\{\bar w^\eps\}_{\eps>0}$ converges in $C\left(K\right)$.
Note first that the solution to problem (\ref{eq:barw}) 
is the sum of the solution to the version of (\ref{eq:barw}) with $b_r(x)\equiv 0$ and  $g(x,t)\equiv 0$  and the solution
to the version of (\ref{eq:barw}) with $F(x,t)\equiv 0$. Since the solution to the latter problem is  continuous in 
$\overline\Pi$, our task  reduces to proving the convergence of $\{\bar w^\eps\}_{\eps>0}$
in  $C\left(K\right)$ in the case when $b_r(x)\equiv 0$ and  $g(x,t)\equiv 0$. 

Let $J^\eps$ be the union, over $i\le n$, $j\le m_i$, and 
$\al\in[-\eps;\eps]$, of all EPs 
emanating from $(x_{ij}^*+\al,0)$.
Fix an arbitrary $i\le n$ and $\eps_0>0$ such that $K\subset\overline\Pi\setminus J^{\eps_0}$.
It is easily seen that, given $(x,t)\in K$,  $\bar w_i^\eps(x,t)$ for 
each $\eps\le\eps_0$ is a finite sum of integrals over the characteristics 
 in the EPs passing through $\om_i(\tau;x,t)$ 
and restricted to $\overline{\Pi^t}\cap J^\eps$
 (see Fig.~4). Denote these by
$C_1^{\eps},\dots,C_{l}^{\eps}$. Moreover, suppose that $C_j^{\eps}$
is given by equation $x=c_j(t)$.

Let $f_{q,s}$ denote the entries of the matrix $F$.
Moreover,  $q(j)$ is the index of the 
characteristic curve passing through $C_j^{\eps}$ and $s(j)$ is 
the  index of the characteristic curves bounding
the ``tube''  containing $C_j^{\eps}$.
 An explicit calculation of $\bar w_i^\eps(x,t)$ gives us equality

\begin{figure}
\centerline{\includegraphics{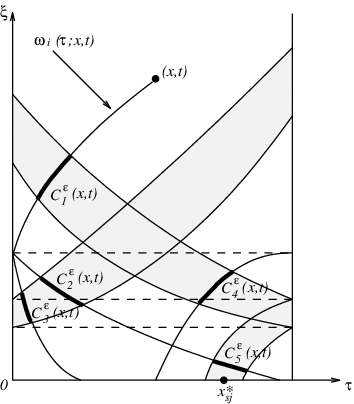}}
\caption{Construction of the set $J_*^\eps$.}
\end{figure}

\begin{equation}\label{eq:Cj}
\bar w_i^\eps(x,t)=\sum\limits_{j=1}^{l}P_{ij}(x,t)\int\limits_{t_j^{\eps-}}^{t_j^{\eps+}}
E_{q(j)}(\tau;c_j(t_j^{\eps+}),t_j^{\eps+})f_{q(j),s(j)}(c_j(\tau),\tau)z_{s(j)}^\eps(c_j(\tau),\tau)\,d\tau,
\end{equation}
where $[t_j^{\eps-},t_j^{\eps+}]$ is the projection of $C_j^{\eps}$ onto $t$-axis and
$P_{ij}(x,t)$ is a smooth function of  $p_{ij}$.

Furthermore,
for any $j\le l$ there are $\al(j)\le n$, $\be(j)\le m_\al$,  $i_1,\dots,i_\ga\in\N$, 
and smooth mappings $T_{i_1} : t\longrightarrow t_{i_1}(c_j(t),t)$, 
$T_{i_p} : t\longrightarrow t_{i_p}(x_{i_{p}},t)$ for $\quad 1< p<\ga$, and
$T_{i_\ga} : t\longrightarrow \om_{i_\ga}(0;x_{i_\ga},t)$,
where $x_i=0$ for $i\le k$ and $x_i=1$ for $k+1\le i\le n$.
Geometrically speaking, for each $j\le l$ there is $x_{\al(j),\be(j)}^*$ and a 
continuous path in $J_*^\eps$ containing both 
the curve $C_{j}^\eps$ and the point $(x_{\al(j),\be(j)}^*,0)$, whatsoever
$\eps\le\eps_0$. 
Moreover, the value of $\ga$  equals to the number of the
``tubes'' creating this path.

Note that $l$ as a function of $(x,t)$ is constant on any connected component of $\overline{\Pi}\setminus J^{\eps_0}$. 
From now on by $l$ we will denote the value of this function on $K$. Calculating of
$z_{s(j)}^\eps$ and changing variables by means of $T_{i_p}$, we transform expression (\ref{eq:Cj})
for $\bar w_i^\eps(x,t)$ on $K$ to a form more convenient for our purpose, namely,
\begin{equation}\label{eq:bar-w-eps}
\bar w_i^\eps(x,t)=\sum\limits_{j=1}^{l}
\int\limits_{x_{\al(j),\be(j)}^*-\eps}^{x_{\al(j),\be(j)}^*+\eps}
Q_{ij}(\xi;x,t)\vphi_{\eps}^{(l_{\al(j),\be(j)})}(\xi-x_{\al(j),\be(j)}^*)\,d\xi,
\end{equation}
where  the kernels $Q_{ij}(\xi;x,t)$ are certain smooth functions. Therefore,
$$
\bar w_i^\eps(x,t)\longrightarrow\sum_{j=1}^{l}
(-1)^{l_{\al(j),\be(j)}}\d_\xi^{(l_{\al(j),\be(j)})}Q_{ij}(x_{\al(j),\be(j)}^*;x,t) \mbox{ as } \eps\to 0,
$$
where the convergence is uniform for all $(x,t)\in K$. 
\end{subproof}
\begin{claim}
$w^\eps\to\bar w+\tilde w$ as $\eps\to 0$ in  $C(K)$, where  $\bar w$ and $\tilde w$
are $C(\overline\Pi\setminus J)$-solutions to problems (\ref{eq:barw}) and (\ref{eq:tildew}),
respectively.
\end{claim}
\begin{subproof}
 On the account of Claim~1,
it is sufficient to prove that $w^\eps-\bar w^\eps-\tilde w\to 0$ as $\eps\to 0$ in  $C(K)$. The difference
$w^\eps-\bar w^\eps-\tilde w$ is a solution to the following problem:
\begin{equation}\label{eq:www}
\begin{array}{rcl}
(\partial_t  + \Lambda\partial_x + A)
\left(w^\eps-\bar w^\eps-\tilde w\right)  &=&  -F(\bar w^\eps-\bar w)\\
\left(w^\eps-\bar w^\eps-\tilde w\right)(x,0)&=& 0\\
\displaystyle
\left(w_i^\eps-\bar w_i^\eps-\tilde w_i\right)(0,t) &=&  \displaystyle\sum\limits_{j=1}^{n}p_{ij}(t)
v_{j}^{w^\eps-\bar w^\eps-\tilde w}(0,t)+r_i\left(t,v^{w^\eps}\right)-
r_i\left(t,v^{\bar w+\tilde w}\right), \\
&&\qquad k+1\le i\le n\\
\displaystyle
\left(w_i^\eps-\bar w_i^\eps-\tilde w_i\right)(1,t) &=&  \displaystyle\sum\limits_{j=1}^{n}p_{ij}(t)
v_{j}^{w^\eps-\bar w^\eps-\tilde w}(0,t)+r_i\left(t,v^{w^\eps}\right)-
r_i\left(t,v^{\bar w+\tilde w}\right), \\
&&\qquad 1\le i\le k.
\end{array}
\end{equation}
Rewrite
\begin{eqnarray}
\lefteqn{
r_i(t,v^{w^\eps})-r_i(t,v^{\bar w+\tilde w})
}\nonumber\\&&
=\int\limits_0^1\left(\nabla_yr_i\right)(t,\sigma v^{w^\eps}+(1-\sigma)v^{\bar w^\eps+\tilde w})\,d\sigma
\cdot v^{w^\eps-\bar w^\eps-\tilde w}
+\left(r_i(t,v^{\bar w^\eps+\tilde w})-r_i(t,v^{\bar w+\tilde w})\right).\nonumber
\end{eqnarray}
Similarly to (\ref{eq:r1}), consider an integral form of (\ref{eq:www}) on $K$ with the boundary term not depending on $v^{w^\eps-\bar w^\eps-\tilde w}$.  
Since $\bar w^\eps\to\bar w$ as $\eps\to 0$ uniformly on any compact subset
of $\overline\Pi\setminus J$,
the difference  $r_i(t,v^{\bar w^\eps+\tilde w}(t))-r_i(t,v^{\bar w+\tilde w}(t))$ tends to zero
as $\eps\to 0$ uniformly on any compact subset of
$(\d\Pi\setminus\{t=0\})\setminus J$.

It remains to show that, for all $i\ne j$, the function 
$$
P_{ij}^\eps(x,t)=\int\limits_{t_i(x,t)}^tE_i(\tau;x,t)
f_{ij}(\xi,\tau)\left(\bar w_j^\eps-\bar w_j\right)(\xi,\tau)|_{\xi=\om_i(\tau;x,t)}\,d\tau
$$
converges to zero as $\eps\to 0$ uniformly on $K$. 
Fix an arbitrary $\eps_1\le\eps_0$ and split up the integral into two parts, namely
$$
P_{ij}^\eps(x,t)=\int\limits_{[t_i(x,t),t]-\triangle_{\eps_1}(x,t)}d\tau+
\int\limits_{\triangle_{\eps_1}(x,t)}d\tau,
$$
where $\triangle_{\eps_1}(x,t)=\{\tau\in[t_i(x,t),t]\,:\,(\om_i(\tau;x,t),\tau)\in J^{\eps_1}\}$.  
Since $\bar w_j^\eps$ is uniformly bounded on each bounded subset of $\overline\Pi$ and
$\bar w\in C(\overline\Pi\setminus J)$, the absolute value of the second integral
is bounded from above by $C\eps_1$ for all $(x,t)\in K$ and for some $C>0$ which does not depend on $\eps_1$.
 Furthermore, there exists $\eps_2\le\eps_1$ such that for all $\eps\le \eps_2$
the absolute value of the first integral is bounded from above by $\eps_1$, 
due to the uniform convergence of $\bar w_j^\eps$ off $J^{\eps_1}$. We
hence conclude that for each $\eps_1\le\eps_0$ there is  $\eps_2\le\eps_1$ and a constant $C>0$ such that for all
$(x,t)\in K$ and $\eps\le\eps_2$ we have $|P_{ij}^\eps(x,t)|\le C\eps_1$. Therefore,
$Q_{ij}^\eps(x,t)\to 0$ as $\eps\to 0$ uniformly on $K$.
\end{subproof}
To finish the proof of this part of the theorem, it remains to recall that problem (\ref{eq:w}) for $w$ has a 
unique solution of the form $\bar w+\tilde w$, where $\bar w$ and $\tilde w$ are the unique solutions to problems
(\ref{eq:barw}) and (\ref{eq:tildew}), respectively.

{\bf 4.}
This part of the theorem is a direct consequence of the three preceding parts.

\end{proof}

In Theorem~\ref{thm:delta} we established the existence of a delta wave solution to 
the problem under consideration. Now we identify conditions under which this solution 
is smoothing. In the following theorem conditions ($\io$) and ($\io\io$) are supposed to be formulated
correspondingly to problem (\ref{eq:1}), (\ref{eq:2'}), (\ref{eq:3'}) .

\begin{theorem}\label{thm:reg-delta}
 Assume that the data $\lambda_i$, $a_{ij}$, $g_i$, and $p_{ij}$ are 
smooth functions in all their arguments and both $r_i$ and $\nabla_yr_i$ are in $\L^\infty\left((0,T)\times\R^n\right)$ 
for any $T>0$. Suppose that condition (\ref{eq:L1}) is fulfilled.
\begin{itemize}
\item
{\rm Sufficiency.} If condition ($\io$) is true, then 
the delta wave solution to problem (\ref{eq:1}), (\ref{eq:2'}), (\ref{eq:3'})  
is smoothing for any  $b_r\in C[0,1]^n$ satisfying equalities  (\ref{eq:zero}).
\item
{\rm Necessity.}
 Assume that  the delta wave solution to problem 
(\ref{eq:1}), (\ref{eq:2'}), (\ref{eq:3'})  
is smoothing for any  $b_r\in C[0,1]^n$ satisfying  equalities  (\ref{eq:zero}).
 Then condition ($\io\io$) is fulfilled.
\end{itemize} 
\end{theorem}

\begin{proof}
{\it Sufficiency.}
As condition ($\io$) is true, there exists 
$T_0$ such that $J_*\in\overline{\Pi^{T_0}}$ and the restriction of $u$ to  $\overline\Pi\setminus\Pi^ {T_0}$
equals $w$.
The smoothingness of $u$ will follow from Theorem~\ref{thm:I} and from the fact that there is $T>T_0$ such that 
$w$ is continuous on $\overline\Pi\setminus\overline\Pi^T$.  
To prove the latter,  fix $T$ to be a supremum over $t>0$ such that there is $x\in[0,1]$ and an EP jointing the point $(x,t)$
with the line $t=T_0$.
Note that $u$ on $\overline\Pi\setminus\overline\Pi^T$ fulfills the system
of integral equations (\ref{eq:u})--(\ref{eq:v}) where
$
h_i(t,v)=\sum_{j=1}^np_{ij}(t)v_j(t)+r_i(t,v).
$
By  the choice of $T$, the right hand side of (\ref{eq:u}) can be rewritten in the form
involving integrals of kind $I_{ijs}(x,t)$ defined by (\ref{eq:I0}). This expression
will not depend neither on $v$ nor on
$u(x,T_0)$. We therefore can use a transformation of $I_{ijs}(x,t)$ similar to that made in
(\ref{eq:I}). Since $u$ is piecewise continuous on $\overline\Pi\setminus\overline\Pi^T$, the right hand side of
(\ref{eq:I}) is a continuous function on $\overline\Pi\setminus\overline\Pi^T$, as desired.

{\it Necessity.}
By Theorem~\ref{thm:delta}, the delta wave solution is given by the sum of a regular part 
$w$ in  $C(\overline\Pi\setminus J)^n$ and a  singular part $z\in\D'(\Pi)$ which, in its turn, is the sum
of strong singularities concentrated on the characteristic curves contributing into the set $J_*$.
By the assumption,  the delta wave solution is smoothing for each $b_r\in C[0,1]^n$. 
This entails the smoothingness of  both  the regular and the singular parts of the solution.
It remains to note that  condition ($\io\io$) is 
necessary  for smoothingness of $w$. This follows by similar argument that was used to prove 
the necessity in Theorem~\ref{thm:I}.
\end{proof}

\section*{Acknowledgments}
This work is supported by a Humboldt Research Fellowship.
The author is thankful to Natal'ya Lyul'ko for helpful discussions.

\end{document}